\documentclass[10pt]{article}
\usepackage{flafter,amsmath,amssymb,latexsym,psfrag,graphicx,color,indentfirst}
\usepackage[margin=2.5cm]{geometry}
\usepackage[numbers,sort&compress]{natbib}

\usepackage[colorlinks,
linktocpage=true,
linkcolor=black,
citecolor=black
]{hyperref}

\newtheorem{theorem}{Theorem}[section]

\newtheorem{lemma}[theorem]{Lemma}
\newtheorem{remark}[theorem]{Remark}

\newtheorem{cor}[theorem]{Corollary}

\newtheorem{example}[theorem]{Example}

\allowdisplaybreaks

\makeatletter

\@addtoreset{figure}{section}
\makeatother

\begin{document}
\setlength\arraycolsep{2pt}
\date{\today}

\title{ Analysis of the acoustic waves reflected by a cluster of small holes in the time-domain and the equivalent mass density}
\author{Mourad Sini$^1$, Haibing Wang$^2$, Qingyun Yao$^2$
\\$^1$RICAM, Austrian Academy of Sciences, A-4040, Linz, Austria \qquad\qquad\quad\quad
\\E-mail: mourad.sini@oeaw.ac.at
\\$^2$School of Mathematics, Southeast University, Nanjing 210096, P.R. China\;
\\E-mail: hbwang@seu.edu.cn
}

\maketitle
\begin{abstract}

We study the time-domain acoustic scattering problem by a cluster of small holes (i.e. sound-soft obstacles). Based on the retarded boundary integral equation method, 
we derive the asymptotic expansion of the scattered field as the size of the holes goes to zero. Under certain geometrical constraints on the size and the minimum distance
of the holes, we show that the scattered field is approximated by a linear combination of point-sources where the weights are given by the capacitance of each hole
and the causal signals (of these point-sources) can be computed by solving a, retarded in time, linear algebraic system. A rigorous justification of the asymptotic expansion and the unique 
solvability of the linear algebraic system are shown under natural conditions on the cluster of holes. As an application of the asymptotic expansion, we derive, in the limit 
case when the holes are densely distributed and occupy a bounded domain, the equivalent effective acoustic medium (an equivalent mass density characterized by the capacitance of the holes) 
that generates, approximately, the same scattered field as the cluster of holes. Conversely, given a locally variable, smooth and positive mass density, satisfying a certain subharmonicity condition, we can 
design a perforated material with holes, having appropriate capacitances, that generates approximately the same acoustic field as the acoustic medium modelled by the given mass density (and constant speed of propagation).
 Finally, we numerically verify the asymptotic expansions by comparing the asymptotic approximations with the 
numerical solutions of the scattered fields via the finite element method.



\end{abstract}

\section{Introduction}\label{intro}
\setcounter{equation}{0}

We are concerned with the time-domain acoustic scattering from a cluster of small sound-soft obstacles (i.e. holes) located in the homogeneous background medium in $\mathbb{R}^3$. 
Let $D$ be a union of holes, i.e. $D=\cup_{j=1}^M D_j$. Assume that $D$ is bounded and $\mathbb R^3\setminus\overline D$ is connected. 
We denote by $c_0$ the constant wave speed in $\mathbb R^3\setminus\overline D$. Let
\begin{equation}\label{ui}
u^i(x,\,t):= \frac{\lambda(t-c_0^{-1}|x-z^*|)}{4\pi |x-z^*|}
\end{equation}
be an incident wave emitted from a point source located at $z^*\not\in \overline D$, where $\lambda\in C^\infty(\mathbb R)$ is a causal signal such that $\lambda$ vanishes for all $t<0$. 
We note that $c_0^{-2} u^i_{tt} - \Delta u^i =0$ for $x\in\mathbb R^3\setminus\{z^*\}$ and $t>0$. Then the scattered acoustic wave $u^s:=u^s(x, t)$, generated after the incident wave hits the cluster $D$, satisfies the following initial boundary value problem:
\begin{equation}\label{ibvp}
\begin{cases}
c_0^{-2}u^s_{tt} - \Delta u^s=0 & \textrm{ in } (\mathbb R^3\setminus\overline D)_T, \\
u^s=-u^i & \textrm{ on } (\partial D)_T, \\
u^s|_{t=0}=0,\;u^s_t|_{t=0}=0 & \textrm{ in } \mathbb R^3\setminus\overline D.
\end{cases}
\end{equation}
 We denote by $u:=u^i+u^s$ the total field. For simplicity of notations, here and throughout this paper, we denote $X\times (0,\,T)$ 
and $\partial X\times (0,\,T)$ by $X_T$ and $(\partial X)_T$, respectively, where $X$ is a domain in $\mathbb R^3$ and $\partial X$ denotes its boundary. The uniqueness and existence of the solution 
to the direct scattering problem \eqref{ibvp} are well studied by using the retarded boundary integral equation method; see for instance \cite{B-H1986_1, Costabel2003, Ha2003, Lubich1994, Sayas2016}.

Now, we consider our holes to be of the form $D_j = \varepsilon B_j + z_j,\,j=1,\,2,\,\cdots,\,M$,  characterized by the parameter $\varepsilon>0$ and the locations $z_j\in \mathbb R^3$, 
where $B_j$'s are bounded and $C^2$-smooth domains containing the origin. The parameter $\varepsilon$ is the relative size of $D_j$, as compared to the size of $B_j$, (i.e. a dimensionless quantity) 
and it is intended to be small. To fix some notations, we set $ a$ as the maximum among the diameters of the holes, i.e.,
\begin{equation}\label{defn_a}
a:=\max_{1\leq j \leq M}\,\mathrm{diam}(D_j)=\varepsilon \, \max_{1\leq j \leq M}\,\mathrm{diam}(B_j),
\end{equation}
and $d$ as the minimum distance between the holes, i.e.,
\begin{equation}\label{defn_d}
d:=\min_{\substack{1\leq i,\,j \leq M\\i\neq j}}\,d_{ij}, \quad d_{ij}:= \mathrm{dist}(D_i,\,D_j).
\end{equation}

Because of \eqref{defn_a}, we sometimes abuse the notations $a$ and $\varepsilon$ when they naturally appear in some estimates. In this paper,
we are interested in the following regimes for modeling the cluster
\begin{equation}\label{M_reg}
M \thicksim \varepsilon^{-s},\; d \thicksim \varepsilon^\beta \quad \mathrm{as}\; \varepsilon \ll 1
\end{equation}
with positive constants $s$ and $\beta$.

\medskip
In this work, we are interested in analyzing the asymptotic behavior of the scattered field $u^s$ for the above time-domain scattering model as the relative size $\varepsilon$ of the holes goes to zero. 
As we know, asymptotic expansions of the fields generated by a cluster of small particles (of different kinds) are well developed in the literature for the elliptic models; 
see for instance \cite{A-C-C-S-18, A-C-C-S-19, B-L-P:1978, C-M-S-17, J-K-O:1994, M-K:2006, A-C-K-S-15, A-K:2007, Bendali-et.al:2015, BS2018, C-S2014, Ramm2015, Martin:2006, M-M-N2017, M-N-P2000, Nie2017} and the references therein. 
However, there are few results on time-domain models, as those related to parabolic, Schr\"odinger or hyperbolic equations, unless for periodic media \cite{C-D2016, D-G2008, D-Z2012, M-K:2006} or finitely many holes \cite{K-P2017}.
Recently, we studied in \cite{S-W2019} the asymptotic analysis of the solution to a heat conduction problem by a cluster of small cavities using the boundary integral equation method, and then derived 
an equivalent effective medium that generate approximately the same temperature field as the cluster of cavities. For the time-domain wave scattering problem, the situation is much less clear. 
This motivates our work in this paper.

Our first result is stated as follows.
\begin{theorem}\label{Main}
We assume that the incident wave is causal and $\lambda \in C^\infty[0,\,T]$ such that $\sum^\infty_{n=0}\mathcal C_n$ is finite where $\mathcal C_n:=\max_{t\in[0,\,T]}\vert \lambda^{(n)}(t)\vert$.  Under the following condition on the cluster of holes:
\begin{equation}\label{M1M_cond_1}
\varepsilon \max_{1\leq i \leq M}\sum_{j\not=i} d_{ij}^{-2} < 1,
\end{equation}
which means that $1-2\beta-s/3\geq 0$, we have the following asymptotic expansion:
\begin{equation}\label{M1M_main}
u^s(x,\,t) = \sum_{j=1}^M \frac{C_j \alpha_j(t-c_0^{-1}|x-z_j|)}{4\pi |x-z_j|} + O\left(\varepsilon^{2-s}\right) + O\left(\varepsilon^{3-2s}\right)+ O\left(\varepsilon^{3-2\beta-s}\right) \quad \mathrm{as}\;\varepsilon\to 0
\end{equation}
for $x\in \mathbb R^3\setminus\overline D$, $x$ away from $D$, and $t\in (0,\,T)$, where the constant $C_j$ is the capacitance of $D_j$ defined by
\begin{equation}\label{C_i}
C_j:=\int_{\partial D_j} \sigma_j(y)\,ds(y), \quad j=1,\,2,\cdots,\,M
\end{equation}
with $\sigma_j$ satisfying
\begin{equation}\label{C_i1}
\int_{\partial D_j} \frac{\sigma_j(y)}{4\pi |x-y|}\,ds(y) = 1 \quad \mathrm{on}\;\partial D_j,
\end{equation}
and $\{\alpha_j \}_{j=1}^M$ is the unique solution of the invertible, retarded in time, linear algebraic system
\begin{equation}\label{M1Mq_i}
\alpha_i(t) + \sum_{\substack{j=1\\j\not= i}}^M \frac{C_j\alpha_j(t-c_0^{-1}|z_i-z_j|)}{4\pi |z_i-z_j|} = - \frac{\lambda(t-c_0^{-1}|z_i-z^*|)}{4\pi |z_i-z^*|}, \quad i=1,\,2,\,\cdots,\,M.
\end{equation}
\end{theorem}

This kind of asymptotic expansions plays a key role in highly important applied sciences, such as imaging and material sciences. 
As an application of the asymptotic expansion in the limit case that the holes are densely distributed and occupy a bounded domain, 
we derive an effective medium that can produce approximately the same scattered field as the union of holes. Explicitly, the solution of the 
initial boundary value problem \eqref{ibvp} can be approximated by the solution of an effective initial value problem whose governing equation contains a zero 
order term generated by the capacitance of the holes. To show this, let $\Omega$ be a bounded domain containing all the holes $D_j,\,j=1,\,2,\,\cdots,\,M$.
Here we assume that the holes have the same shapes\footnote{Actually, we only need them to have the same capacitance.}. We know that the capacitance $C_j$
 of $D_j$ is given by the one of $B_j$, which we denote by $\overline{C}_j$, through the formula $C_j = \overline{C}_j \varepsilon$. We set the common capacitance $\overline{C}_j$ by $\overline{C}$.
 We divide $\Omega$ into $[a^{-1}]$ subdomains $\Omega_j,\,j=1,\,2,\,\cdots,\,[a^{-1}]$, periodically arranged for instance\footnote{The periodicity is actually not needed.}, such that the $\Omega_j$'s are disjoint 
and of a volume $a$. Let each subdomain $\Omega_j$ contain one hole. Such a distribution obeys the condition (\ref{M1M_cond_1}), with $\beta=\frac{1}{3}$ and $s=1$, as we explained in \cite{S-W2019}. 
Then we have the following result.

\begin{theorem}\label{Main2}
Let $W(x,\,t)$ be the solution of the initial value problem
\begin{equation}\label{MW-1}
\begin{cases}
(c_0^{-2}\partial_{tt} - \Delta + \overline C \chi_\Omega)W = - \overline C \chi_\Omega u^i(x,\,t) & \mathrm{in}\;\mathbb R^3\times (0,\,T), \\
W(x,\,0)=0,\,W_t(x,\,0)=0 & \mathrm{in}\; \mathbb R^3.
\end{cases}
\end{equation}
Then, for any fixed $x\in\mathbb R^3\setminus\overline \Omega$ and $t\in (0,\,T)$, we have the estimate
\begin{equation}\label{Meff}
u^s(x,\,t) = W(x,\,t) + O\left( \varepsilon^{\frac{1}{3}} \right)  \quad \mathrm{as}\; \varepsilon\to 0,
\end{equation}
where $u^s(x,\,t)$ is the solution to \eqref{ibvp}. If we define $U:=W+u^i$, we also have
\begin{equation}\label{Meff1}
u(x,\,t) = U(x,\,t) + O\left( \varepsilon^{\frac{1}{3}} \right)  \quad \mathrm{as}\;\varepsilon\to 0.
\end{equation}
\end{theorem}

Let $p$ be the unique solution of the problem
\begin{equation}\label{p}
\begin{cases}
-\Delta p +\overline C p =0 & \mathrm{ in }\; \Omega,\\
p =1 & \mathrm{ on }\; \partial \Omega.
\end{cases}
\end{equation}
Since $\overline C$ is positive in $\Omega$, the unique solution of \eqref{p} is also positive in $\Omega$, due to the maximum principle. We extend $p$ from $\Omega$ to $\mathbb{R}^3$ by simply setting $p=1$ in $\mathbb{R}^3\setminus{\overline \Omega}$. Define $\tilde{U}:=p^{-1}U$ and observe that $\tilde U$ satisfies the problem
\begin{equation}\label{tilde-U}
\begin{cases}
p^2 c_0^{-2}\partial_{tt} \tilde U - \nabla \cdot \left(p^2 \nabla \tilde U\right) = p\, \lambda(t)\, \delta(x-z^*) & \mathrm{in}\;\mathbb R^3\times (0,\,T), \\
\tilde U(x,\,0)=0,\; \tilde U_t(x,\,0)=0 & \mathrm{in}\; \mathbb R^3.
\end{cases}
\end{equation}
Then, as a corollary of Theorem \ref{Main2}, we deduce the following result.
\begin{cor}
For $x\in\mathbb R^3\setminus\overline \Omega$ and $t\in (0,\,T)$, we have the approximation
\begin{equation}\label{Meff2}
u(x,\,t) = \tilde U(x,\,t) + O\left( \varepsilon^{\frac{1}{3}} \right)  \quad \mathrm{as}\,\varepsilon\to 0.
\end{equation}
\end{cor}
This result means that the wave reflected by the cluster of holes is approximately the same as the one generated by the acoustic medium characterized by the speed of propagation $c_0$ 
and the mass density $\rho:=p^{-2}$, where $p$ is the unique solution of the problem (\ref{p}). Conversely, let $\rho$ be any given mass density function which is $C^2$-smooth, positive 
such that $p:=\rho^{-1/2}$ is subharmonic, i.e. $\Delta p>0$ in a given region $\Omega$ and $p=1$ in $\mathbb{R}^3\setminus\Omega$. Then starting from a homogeneous material (i.e. the background), we drill small holes 
$D_j$ of center $z_j$ and radius $\varepsilon$ having the capacitance $\overline{C}(z_j)\; \varepsilon$, where $\overline{C}:=\frac{\Delta p}{p}$, distributed, periodically for instance, in $\Omega$. 
This perforated material will behave as an acoustic medium with constant speed of propagation and the mass density as the given function $\rho$.\footnote{The proof of Theorem \ref{Main2} is proved for holes having the same capacitances.
However, we do believe that the same result is true for variable capacitances as described above; see \cite{A-C-K-S-15} for the time harmonic acoustic model.}

Here we would like to add the following two observations:
\begin{enumerate}
 
\item The periodicity in distributing the holes in $\Omega$ is actually not needed. We assume it only for simplicity of exposition, and the result can be extended to more general cases. 
In addition, we can put arbitrary number of holes in each subdomain $\Omega_j$. In this case, we need to introduce the local distribution density function $K(x)$ and replace
$\overline C$ by $\overline C\,K(x)$ in the governing equation; see for instance \cite{A-C-K-S-15, A-C-C-S-19, C-M-S-17} for the harmonic regime cases.
    
\item Theorem \ref{Main2} may have important applications in material sciences. On one hand, we can design new materials by appropriately distributing 
the holes in the background medium so that we can get the desired mass density and the scattered field as we explained above. On the other hand, 
for the wave scattering from an inhomogeneous medium modeled by $(c_0^{-2}\partial_{tt} - \Delta - q(x))u^s=0$, we may kill the term $q(x)u^s$ and make the 
inhomogeneity invisible through properly embedding the holes into the medium such that $\overline C\,K(x)=q(x)$. This would be a good insight for acoustic cloaking in the time-domain.
    
\end{enumerate}

The rest of the paper is organized as follows. In Section \ref{single}, we provide the analysis for the case of a single hole to describe the main steps of our approach. 
In Section \ref{multiple}, we prove the asymptotic expansion for the case of multiple holes, i.e. Theorem \ref{Main}. In Section \ref{effective-medium}, 
we derive the effective medium and prove Theorem \ref{Main2}. Three numerical examples are presented in Section \ref{numer} to illustrate the effectiveness of the asymptotic expansion.

\section{Proof of Theorem \ref{Main}: the single hole case}\label{single}
\setcounter{equation}{0}

In this section, we consider the single hole case that $D=\varepsilon\, B + z$ and prove Theorem \ref{Main} (with $M=1$). To begin with, we introduce the function space
\begin{equation*}
H_0^r(0,\,T):=\left\{ g|_{(0,\,T)}:\; g\in H^r(\mathbb R)\, \textrm{ with }\, g\equiv 0 \, \textrm{ in }\, (-\infty,\,0) \right\},\quad r\in\mathbb R
\end{equation*}
and generalize it to the $H^s(\partial D)$-valued function space, denoted by $H_0^r(0,\,T;\, H^s(\partial D))$. Let $E$ be a Hilbert space and define
\begin{equation*}
LT(\sigma,\,E):= \left\{f\in \mathcal D^\prime_+(E):\,e^{-\sigma t}f \in \mathcal S^\prime_+(E) \right\}, \quad \sigma>0,
\end{equation*}
where $\mathcal D^\prime_+(E)$ and $\mathcal S^\prime_+(E)$ denote the sets of distributions and temperate distributions on $\mathbb R$ with values in $E$ and support in $[0,\,+\infty)$. Then we define the space
\begin{equation*}
H_{0,\sigma}^r(0,\,T;\, H^s(\partial D)):= \left\{f\in LT(\sigma,\,H^s(\partial D)) : \, e^{-\sigma t}\Lambda^r f\in L^2_0(0,\,T;\,H^s(\partial D))  \right\},
\end{equation*}
where $r\in\mathbb R$ and $\Lambda^r$ denotes the $r$-th order derivative with respect to the variable $t$. For nonnegative integer $r$, we use the norm
\begin{equation*}
\| f \|_{H_{0,\sigma}^r(0,\,T;\, H^s(\partial D))}:=\left(\int_0^T  e^{-2\sigma t}\left[\left\| f \right\|^2_{H^s(\partial D)} + \sum_{k=1}^r \left\| \frac{\partial^k f}{\partial t^k} \right\|^2_{H^s(\partial D)}\right] \,dt\right)^{1/2}.
\end{equation*}

We now express the solution to \eqref{ibvp} as a retarded single-layer potential
\begin{equation}\label{layer}
u^s(x,\,t)=\int_{\partial D} \frac{\varphi(y,\,t-c_0^{-1}|x-y|)}{4\pi |x-y|}\,ds(y), \quad (x,\,t)\in(\mathbb R^3\setminus\overline D)_T,
\end{equation}
where $\varphi$ is a causal density to be determined. In view of the boundary condition in \eqref{ibvp}, we obtain from the continuity property of the potential \eqref{layer} that
\begin{equation}\label{bie}
\int_{\partial D} \frac{\varphi(y,\,t-c_0^{-1}|x-y|)}{4\pi |x-y|}\,ds(y) = -u^i(x,\,t), \quad (x,\,t)\in (\partial D)_T.
\end{equation}
It was proved in \cite{Lubich1994} that the boundary integral equation \eqref{bie} has a unique solution with the {\it a-priori} estimate
\begin{equation}\label{est_sigma_1}
\| \varphi \|_{H_{0,\sigma}^r(0,\,T;\,H^{-1/2}(\partial D))} \lesssim \| u^i \|_{H_{0,\sigma}^{r+2}(0,\,T;\,H^{1/2}(\partial D))},\quad r\in\mathbb R.
\end{equation}
By the embedding $H^r(0,\,T)\hookrightarrow C[0,\,T]$ for $r>1/2$, we also have
\begin{equation}\label{est_sigma_2}
\| \varphi(\cdot,\,t) \|_{H^{-1/2}(\partial D)} \lesssim \| u^i \|_{H_{0,\sigma}^r(0,\,T;\,H^{1/2}(\partial D))}, \quad r>5/2,\;t\in [0,\,T].
\end{equation}
Throughout the paper, we use the notation \lq\lq$\lesssim$\rq\rq{} to denote \lq\lq$\leq$\rq\rq{} with its right-hand side multiplied by a generic positive constant, if we do not emphasize the dependence of the constant on some parameters.

As we need to deal with changes of coordinates in estimating $\varphi$ by Sobolev norms, we introduce some notations here. We first consider the scaling for the space variable. Set
\begin{equation*}
\hat p(\xi)=p^\wedge(\xi):=p(\varepsilon\xi+z), \; \xi \in \partial B \quad \textrm{and} \quad
\check q(x)=q^\vee(x):=q\left(\frac{x-z}{\varepsilon}\right), \; x\in \partial D.
\end{equation*}
We introduce the following Sobolev norms defined in \cite{G-R1986}:
\begin{equation}\label{defn_Hnorm1}
\|p\|_{H^{1/2}(\partial D)}:=\inf_{\substack{u\in H^1(D)\\u|_{\partial D}=p}}\| u\|_{H^1(D)} \quad \textrm{ for all } p\in H^{1/2}(\partial D)
\end{equation}
and
\begin{equation}\label{defn_Hnorm2}
\|q\|_{H^{-1/2}(\partial D)}:= \sup_{\substack{p\in H^{1/2}(\partial D)\\ p\neq 0}} \frac{|\langle q,\,p \rangle_{\partial D}|}{\|p\|_{H^{1/2}(\partial D)}} \quad \textrm{ for all } q\in H^{-1/2}(\partial D),
\end{equation}
where $\langle \cdot,\,\cdot \rangle_{\partial D}$ denotes the duality paring between $H^{-1/2}(\partial D)$ and $H^{1/2}(\partial D)$. Let
\begin{equation*}
H^{1/2}_\diamond(\partial D):=\big\{p\in H^{1/2}(\partial D):\,\int_{\partial D} p\,ds =0 \big\}\quad \textrm{and} \quad
H^{-1/2}_\diamond(\partial D):=\big\{q\in H^{-1/2}(\partial D):\,\int_{\partial D} q\,ds =0 \big\}.
\end{equation*}
Then we have the following properties for scaling the space variable.

\begin{lemma}\label{scaling}
Suppose $0< \varepsilon \leq 1$.
If $p \in H^{1/2}_\diamond(\partial D)$ and $q \in H^{-1/2}_\diamond(\partial D)$, there exist two constants $c_1$ and $c_2$ such that
\begin{eqnarray}
&&c_1\,\varepsilon^{1/2} \|\hat p\|_{H^{1/2}(\partial B)} \leq \|p\|_{H^{1/2}(\partial D)} \leq  \varepsilon^{1/2} \|\hat p\|_{H^{1/2}(\partial B)}, \label{scaling_11} \\
&&\varepsilon^{3/2} \|\hat q\|_{H^{-1/2}(\partial B)} \leq \|q\|_{H^{-1/2}(\partial D)} \leq c_2\, \varepsilon^{3/2} \|\hat q\|_{H^{-1/2}(\partial B)}. \label{scaling_12}
\end{eqnarray}
If $p\in H^{1/2}(\partial D)$ and $q\in H^{-1/2}(\partial D)$ are constants, there exist two constants $c_3$ and $c_4$ such that
\begin{eqnarray}
&&c_3\, \varepsilon^{3/2} \|\hat p\|_{H^{1/2}(\partial B)} \leq \|p\|_{H^{1/2}(\partial D)} \leq  \varepsilon^{3/2} \|\hat p\|_{H^{1/2}(\partial B)}, \label{scaling_21}\\
&&\varepsilon^{1/2} \|\hat q\|_{H^{-1/2}(\partial B)} \leq \|q\|_{H^{-1/2}(\partial D)} \leq  c_4\, \varepsilon^{1/2} \|\hat q\|_{H^{-1/2}(\partial B)}. \label{scaling_22}
\end{eqnarray}
\end{lemma}
{\bf Proof.} The scaling results \eqref{scaling_11} and \eqref{scaling_12} were proved in \cite[Lemma 4.1]{A-G-H2007}, 
while \eqref{scaling_21} and \eqref{scaling_22} can also be observed from the proof there. \hfill $\Box$

\medskip
Next, we do the scaling for both the space and time variables. Denote $T_\varepsilon:= T/\varepsilon$. For any functions $\varphi$ and $\psi$ defined on $(\partial D)_T$ 
and $(\partial B)_{T_\varepsilon}$, respectively, we use the notations
\begin{eqnarray*}
&&\hat\varphi(\xi,\,\tau)=\varphi^\wedge(\xi,\,\tau):=\varphi(\varepsilon\xi+z,\,\varepsilon \tau), \quad (\xi,\,\tau)\in (\partial B)_{T_\varepsilon},\\
&&\check\psi(x,\,t)=\psi^\vee(x,\,t):=\psi\left(\frac{x-z}{\varepsilon},\,\frac{t}{\varepsilon}\right), \quad (x,\,t)\in (\partial D)_T.
\end{eqnarray*}
Notice that
\begin{equation*}
\frac{\partial^n \hat\varphi(\cdot,\,\tau)}{\partial \tau^n}=\frac{\partial^n \varphi(\cdot,\,\varepsilon \tau)}{\partial \tau^n} = \varepsilon^n \frac{\partial^n \varphi(\cdot,\,t)}{\partial t^n}, \quad n \in\mathbb Z_+.
\end{equation*}
Then, using Lemma \ref{scaling}, we have the following scaling result.

\begin{lemma}\label{scalingt}
Suppose $0< \varepsilon \leq 1$. If $\psi \in H^{r+2}_{0,\sigma}(0,\,T;\, H^{1/2}_\diamond(\partial D))$ and $\varphi \in H^r_{0,\sigma}(0,\,T;\, H^{-1/2}_\diamond(\partial D))$ with nonnegative integer $r$, there exist two constants $c_1$ and $c_2$ such that
\begin{eqnarray}
&&c_1\,\varepsilon^{-(r+1)} \|\hat \psi\|_{H^{r+2}_{0,\varepsilon\sigma}(0,\,T_\varepsilon;\, H^{1/2}(\partial B))} \leq \|\psi\|_{H^{r+2}_{0,\sigma}(0,\,T;\, H^{1/2}(\partial D))} \leq C(T)\, \varepsilon^{-(r+1)} \|\hat \psi\|_{H^{r+2}_{0,\varepsilon\sigma}(0,\,T_\varepsilon;\, H^{1/2}(\partial B))}, \label{scalingt_11} \\
&&\varepsilon^{2-r} \|\hat \varphi\|_{H^r_{0,\varepsilon\sigma}(0,\,T_\varepsilon;\,H^{-1/2}(\partial B))} \leq \|\varphi\|_{H^r_{0,\sigma}(0,\,T;\,H^{-1/2}(\partial D))} \leq c_2 \,C(T)\, \varepsilon^{2-r} \, \|\hat \varphi\|_{H^r_{0,\varepsilon\sigma}(0,\,T_\varepsilon;\,H^{-1/2}(\partial B))}, \label{scalingt_12}
\end{eqnarray}
where $C(T)$ stands for a constant that depends on $T$. If $\psi \in H^{r+2}_{0,\sigma}(0,\,T;\, H^{1/2}(\partial D))$ and $\varphi \in H^r_{0,\sigma}(0,\,T;\, H^{-1/2}(\partial D))$ are independent of the space variable, there exist two constants $c_3$ and $c_4$ such that
\begin{eqnarray}
&&c_3\, \varepsilon^{-r} \|\hat \psi\|_{H^{r+2}_{0,\varepsilon\sigma}(0,\,T_\varepsilon;\, H^{1/2}(\partial B))} \leq \|\psi\|_{H^{r+2}_{0,\sigma}(0,\,T;\, H^{1/2}(\partial D))} \leq C(T)\, \varepsilon^{-r} \|\hat \psi\|_{H^{r+2}_{0,\varepsilon\sigma}(0,\,T_\varepsilon;\, H^{1/2}(\partial B))}, \label{scalingt_21}\\
&&\varepsilon^{1-r} \|\hat \varphi\|_{H^r_{0,\varepsilon\sigma}(0,\,T_\varepsilon;\,H^{-1/2}(\partial B))} \leq \|\varphi\|_{H^r_{0,\sigma}(0,\,T;\,H^{-1/2}(\partial D))} \leq c_4\,C(T)\,\varepsilon^{1-r}\, \|\hat \varphi\|_{H^r_{0,\varepsilon\sigma}(0,\,T_\varepsilon;\,H^{-1/2}(\partial B))}.  \label{scalingt_22}
\end{eqnarray}
\end{lemma}

{\bf Proof.} We only prove \eqref{scalingt_11}, since the others can be proved in the same way. Note that $H_{0,\sigma}^{r+2}$-norm with respect to $t$ is equivalent to $L^2$-norm of the highest derivative. Then, for $\psi\in H^{r+2}_{0,\sigma}(0,\,T;\,H^{1/2}_\diamond(\partial D))$, we derive
\begin{eqnarray*}
\|\psi\|^2_{H^{r+2}_{0,\sigma}(0,\,T;\,H^{1/2}(\partial D))} &\leq & C(T) \int_0^T e^{-2\sigma t}\, \left\|\frac{\partial^{r+2} \psi(\cdot,\,t)}{\partial t^{r+2}}\right\|^2_{H^{1/2}(\partial D)}\,dt\\
& \lesssim & C(T) \int_0^{T_\varepsilon} e^{-2\sigma \varepsilon \tau}\, \varepsilon^{-2r-3}\, \left\|\frac{\partial^{r+2} \hat\psi(\cdot,\,\tau)}{\partial \tau^{r+2}}\right\|^2_{H^{1/2}(\partial B)} \,\varepsilon d\tau \\
&\lesssim & C(T)\, \varepsilon^{-2(r+1)} \|\hat \psi\|^2_{H^{r+2}_{0,\varepsilon\sigma}(0,\,T_\varepsilon;\,H^{1/2}(\partial B))}.
\end{eqnarray*}
On the other hand, we have
\begin{eqnarray*}
\|\psi\|^2_{H^{r+2}_{0,\sigma}(0,\,T;\,H^{1/2}(\partial D))} &\geq& \int_0^T e^{-2\sigma t}\, \left\|\frac{\partial^{r+2} \psi(\cdot,\,t)}{\partial t^{r+2}}\right\|^2_{H^{1/2}(\partial D)}\,dt\\
&\gtrsim &\int_0^{T_\varepsilon} e^{-2\sigma \varepsilon \tau}\, \varepsilon^{-2r-3} \, \left\|\frac{\partial^{r+2} \hat\psi(\cdot,\,\tau)}{\partial \tau^{r+2}}\right\|^2_{H^{1/2}(\partial B)} \,\varepsilon d\tau \\
&\gtrsim & \varepsilon^{-2(r+1)} \|\hat \psi\|^2_{H^{r+2}_{0,\varepsilon\sigma}(0,\,T_\varepsilon;\,H^{1/2}(\partial B))}.
\end{eqnarray*}
The proof is complete. \hfill $\Box$

\medskip
In the next lemma, we investigate scaling property of the retarded single-layer potential operator $\mathcal S_{\partial D}$ defined by
\begin{equation}\label{S}
\mathcal S_{\partial D}[\varphi](x,\,t):=\int_{\partial D} \frac{\varphi(y,\,t-c_0^{-1}|x-y|)}{4\pi |x-y|}\,ds(y), \quad (x,\,t)\in (\partial D)_T.
\end{equation}
For this purpose, we also define
\begin{equation}\label{SB}
\mathcal S^\varepsilon_{\partial B}[\psi](\xi,\,\tau):=\int_{\partial B} \frac{\psi(\eta,\,\tau-c_0^{-1}|\xi-\eta|)}{4\pi |\xi-\eta|}\,ds(\eta), \quad (\xi,\,\tau)\in (\partial B)_{T_\varepsilon}.
\end{equation}

\begin{lemma}\label{inv_SB}
The inverse of the operator $\mathcal S_{\partial B}^\varepsilon:\,H^r_{0,\varepsilon\sigma}(0,\,T_\varepsilon;\,H^{-1/2}(\partial B)) \to H^{r+2}_{0,\varepsilon\sigma}(0,\,T_\varepsilon;\,H^{1/2}(\partial B))$ with $r=0,\,1$ is estimated by $O(\frac{1}{\varepsilon^{r+1}})$ for $\varepsilon\ll 1$.
\end{lemma}

{\bf Proof.} Denote by $V(s)$ the single-layer potential operator for the Helmholtz equation $\Delta U - s^2\,U=0$. Then, by Proposition 3 in \cite{B-H1986_1}, the operator $V(s):\,H^{-1/2}(\partial B) \to H^{1/2}(\partial B)$ is an isomorphism and the operator norm of its inverse $V^{-1}(s)$ is bounded by
\begin{equation}\label{SB_1}
\| V^{-1}(s) \| \leq  \frac{|s|^2}{c\min\{1,\,\sigma\}}, \quad \mathrm{Re}\,s=\sigma>0,
\end{equation}
where $c$ is a positive constant which depends only on $\partial B$. View the single-layer potential operator $\mathcal S_{\partial B}^\varepsilon$ as a convolution with respect to the time variable, and use the operational notation $K(\partial_t)g:=k*g$, where $k$ is the inverse Laplace transform of $K$ and $g$ vanishes for $t<0$. Then we have $\mathcal S_{\partial B}^\varepsilon = V(\partial_t)$. By Lemma 2.1 in \cite{Lubich1994}, we obtain from \eqref{SB_1} that $V^{-1}(\partial_t)$ extends to a bounded linear operator from $H_{0,\varepsilon\sigma}^{r+2}(0,\,T_\varepsilon)$ into $H_{0,\varepsilon\sigma}^r(0,\,T_\varepsilon)$ for arbitrary real number $r$. More explicitly, for $r=0$ and $g\in H^r_{0,\varepsilon\sigma}(0,\,T_\varepsilon;\,H^{-1/2}(\partial B))$ compactly supported with respect to $t$ in $[0,\,T_\varepsilon]$, we have
\begin{eqnarray*}\label{SB_2}
\left\| \left( \mathcal S^\varepsilon_{\partial B} \right)^{-1} [g]\right\|^2_{H^0_{0,\varepsilon\sigma}(0,\,T_\varepsilon;\,H^{-1/2}(\partial B))} &=&\int_0^{T_\varepsilon} \left(e^{-\varepsilon\sigma t} \,\|V^{-1}(\partial_t)g \|_{H^{-1/2}(\partial B)}\right)^2\,dt\nonumber \\
& \leq & \int_0^{+\infty}\left(\mathcal F \left[e^{-\varepsilon\sigma t} \|V^{-1}(\partial_t)g \|_{H^{-1/2}(\partial B)}\right](\eta)\right)^2\,d\eta \nonumber\\
& = & \int_{\varepsilon\sigma-i\mathbb R_+}\left(\mathcal L \left[ \|V^{-1}(\partial_t)g \|_{H^{-1/2}(\partial B)}\right](s)\right)^2\,ds \nonumber\\
& = & \int_{\varepsilon\sigma-i\mathbb R_+}\|\mathcal L \left[ V^{-1}(\partial_t)g\right](s) \|^2_{H^{-1/2}(\partial B)}\,ds \nonumber\\
& = & \int_{\varepsilon\sigma-i\mathbb R_+}\| V^{-1}(s)\,\mathcal L[g](s) \|^2_{H^{-1/2}(\partial B)}\,ds \nonumber\\
& \lesssim & \frac{1}{(\varepsilon\sigma)^2}\,\int_{\varepsilon\sigma-i\mathbb R_+}\| s^2 \mathcal L[g](s) \|^2_{H^{1/2}(\partial B)}\,ds \nonumber\\
& = & \frac{1}{(\varepsilon\sigma)^2}\,\int_{\varepsilon\sigma-i\mathbb R_+} \| \mathcal L[\partial_t^2 g](s) \|^2_{H^{1/2}(\partial B)}\,ds \nonumber\\
& = & \frac{1}{(\varepsilon\sigma)^2}\,\int_0^{+\infty} \left( e^{-\varepsilon\sigma t}\,\| \partial_t^2 g \|_{H^{1/2}(\partial B)}\right)^2\,dt \nonumber\\
& \leq & \frac{1}{(\varepsilon\sigma)^2}\, \left\| g\right\|^2_{H^2_{0,\varepsilon\sigma}(0,\,T_\varepsilon;\,H^{1/2}(\partial B))} .
\end{eqnarray*}
Similarly, we derive for $r=1$ that
\begin{eqnarray*}\label{SB_3}
\left\| \left( \mathcal S^\varepsilon_{\partial B} \right)^{-1} [g]\right\|^2_{H^1_{0,\varepsilon\sigma}(0,\,T_\varepsilon;\,H^{-1/2}(\partial B))}&=& \int_0^{T_\varepsilon}e^{-2\varepsilon\sigma t} \left( \| V^{-1}(\partial_t)g \|^2_{H^{-1/2}(\partial B)} + \left\|\frac{\partial V^{-1}(\partial_t)g}{\partial t}\right\|^2_{H^{-1/2}(\partial B)}\right)\,dt \nonumber\\
&\leq & (1+T_\varepsilon^2)\int_0^{T_\varepsilon}e^{-2\varepsilon\sigma t}\, \|\partial_t\left( V^{-1}(\partial_t)g\right)\|^2_{H^{-1/2}(\partial B)}  \,dt\nonumber\\
& \leq & (1+T_\varepsilon^2) \int_0^{+\infty} \left(\mathcal F \left[e^{-\varepsilon\sigma t}\, \| \partial_t\left(V^{-1}(\partial_t)g\right)\|_{H^{-1/2}(\partial B)} \right](\eta)\right)^2 \,d\eta \nonumber\\
& \leq & (1+T_\varepsilon^2) \frac{1}{(\varepsilon\sigma)^2}\,\int_{\varepsilon\sigma-i\mathbb R_+} \|s^3 \mathcal L[g](s) \|^2_{H^{1/2}(\partial B)}\,ds \nonumber\\
& = & (1+T_\varepsilon^2) \frac{1}{(\varepsilon\sigma)^2}\,\int_{\varepsilon\sigma-i\mathbb R_+}\| \mathcal L[\partial_t^3 g](s) \|^2_{H^{1/2}(\partial B)}\,ds \nonumber\\
& = & (1+T_\varepsilon^2) \frac{1}{(\varepsilon\sigma)^2}\,\int_0^{+\infty} e^{-2\varepsilon\sigma t}\, \| \partial_t^3 g \|^2_{H^{1/2}(\partial B)}\,dt\nonumber \\
&\leq & (1+T_\varepsilon^2) \frac{1}{(\varepsilon\sigma)^2}\, \left\|g\right\|^2_{H^3_{0,\varepsilon\sigma}(0,\,T_\varepsilon;\,H^{1/2}(\partial B))}.
\end{eqnarray*}
The proof is complete. \hfill $\Box$

\begin{lemma}\label{scalingS}
Let $\varphi\in H^r_{0,\sigma}(0,\,T;\,H^{-1/2}(\partial D))$ and $\psi\in H^{r+2}_{0,\sigma}(0,\,T;\,H^{1/2}(\partial D))$. Then
\begin{eqnarray}
\mathcal S_{\partial D}[\varphi] &=& \varepsilon (\mathcal S^\varepsilon_{\partial B}[\hat\varphi])^\vee, \label{scalingS_1}\\
(\mathcal S_{\partial D})^{-1}[\psi] &=& \varepsilon^{-1} \left((\mathcal S^\varepsilon_{\partial B})^{-1}[\hat\psi]\right)^\vee, \label{scalingS_2}
\end{eqnarray}
and
\begin{eqnarray}\label{scalingS_3}
&&\left \| (\mathcal S_{\partial D})^{-1} \right\|_{\mathcal L \left(H^{r+2}_{0,\sigma}(0,\,T;\,H^{1/2}_\diamond(\partial D)),\,H^r_{0,\sigma}(0,\,T;\,H^{-1/2}(\partial D))\right)} \nonumber\\
&\lesssim& \varepsilon\, \left \| (\mathcal S^\varepsilon_{\partial B})^{-1} \right\|_{\mathcal L \left(H^{r+2}_{0,\varepsilon\sigma}(0,\,T_\varepsilon;\,H^{1/2}_\diamond(\partial B)),\,H^r_{0,\varepsilon\sigma}(0,\,T_\varepsilon;\,H^{-1/2}(\partial B))\right)}.
\end{eqnarray}
\end{lemma}

{\bf Proof.} Let $x=\varepsilon \xi + z$, $y=\varepsilon \eta + z$ and $t=\varepsilon \tau$. Then we have
\begin{eqnarray*}
\mathcal S_{\partial D}[\varphi](x,\,t) &=& \int_{\partial D} \frac{\varphi(y,\,t-c_0^{-1}|x-y|)}{4\pi |x-y|}\,ds(y)\\
&=& \varepsilon  \int_{\partial B} \frac{\varphi(\varepsilon \eta+z,\,\varepsilon \tau-c_0^{-1}\varepsilon|\xi-\eta|)}{4\pi |\xi-\eta|}\,ds(\eta)\\
&=& \varepsilon \mathcal S^\varepsilon_{\partial B}[\hat \varphi](\xi,\,\tau),
\end{eqnarray*}
which gives \eqref{scalingS_1}. Further, the identity \eqref{scalingS_2} follows from the derivation
\begin{equation*}
\mathcal S_{\partial D}\left[\left((\mathcal S^\varepsilon_{\partial B})^{-1}[\hat\psi]\right)^\vee\right] = \varepsilon \left(\mathcal S^\varepsilon_{\partial B}\left[(\mathcal S^\varepsilon_{\partial B})^{-1}[\hat\psi]\right]\right)^\vee = \varepsilon (\hat \psi)^\vee = \varepsilon \psi.
\end{equation*}

To show the estimate \eqref{scalingS_3}, we derive that
\begin{eqnarray*}
& &\left \| (\mathcal S_{\partial D})^{-1} \right\|_{\mathcal L \left(H^{r+2}_{0,\sigma}(0,\,T;\,H^{1/2}_\diamond(\partial D)),\,H^r_{0,\sigma}(0,\,T;\,H^{-1/2}(\partial D))\right)}\\
&:=& \sup\limits_{0\not=\psi\in H^{r+2}_{0,\sigma}(0,\,T;\,H^{1/2}_\diamond(\partial D))} \displaystyle \frac{ \| (\mathcal S_{\partial D})^{-1}[\psi] \|_{H^r_{0,\sigma}(0,\,T;\,H^{-1/2}(\partial D))}}{\|\psi\|_{H^{r+2}_{0,\sigma}(0,\,T;\,H^{1/2}(\partial D))}}\\
&\lesssim & \sup\limits_{0\not=\psi\in H^{r+2}_{0,\sigma}(0,\,T;\,H^{1/2}_\diamond(\partial D))} \displaystyle \frac{ \varepsilon^{1-r} \| \left((\mathcal S_{\partial D})^{-1}[\psi]\right)^\wedge \|_{H^r_{0,\varepsilon\sigma}(0,\,T_\varepsilon;\,H^{-1/2}(\partial B))}}{\varepsilon^{-(r+1)} \|\hat\psi\|_{H^{r+2}_{0,\varepsilon\sigma}(0,\,T_\varepsilon;\,H^{1/2}(\partial B))}}\\
&\lesssim & \sup\limits_{0\not=\hat\psi\in H^{r+2}_{0,\varepsilon\sigma}(0,\,T_\varepsilon;\,H^{1/2}_\diamond(\partial B))}  \displaystyle \frac{ \varepsilon \| (\mathcal S^\varepsilon_{\partial B})^{-1}[\hat\psi] \|_{H^r_{0,\varepsilon\sigma}(0,\,T_\varepsilon;\,H^{-1/2}(\partial B))}}{\|\hat\psi\|_{H^{r+2}_{0,\varepsilon\sigma}(0,\,T_\varepsilon;\,H^{1/2}(\partial B))}}\\
&=& \varepsilon \left \| (\mathcal S^\varepsilon_{\partial B})^{-1} \right\|_{\mathcal L \left(H^{r+2}_{0,\varepsilon\sigma}(0,\,T_\varepsilon;\,H^{1/2}_\diamond(\partial B)),\,H^r_{0,\varepsilon\sigma}(0,\,T_\varepsilon;\,H^{-1/2}(\partial B))\right)}.
\end{eqnarray*}
Thus, the proof is complete. \hfill $\Box$

Here we point out that if we restrict $(\mathcal S_{\partial D})^{-1}$ into the subset consisting of functions independent of 
the space variable in $H^{r+2}_{0,\sigma}(0,\,T;\,H^{1/2}(\partial D))$, then the estimate \eqref{scalingS_3} should be 
$\left \| (\mathcal S_{\partial D})^{-1} \right\| \lesssim \left \| (\mathcal S^\varepsilon_{\partial B})^{-1} \right\|$.

\medskip
We now show an {\it a-priori} estimate of the solution $\varphi$ to \eqref{bie}. Set
\begin{equation*}
g_1(t) := \frac{-1}{|\partial D|} \int_{\partial D} u^i(x,\,t)\,ds(x)\quad \textrm{and} \quad g_2(x,\,t):= -u^i(x,\,t) - g_1(t), \quad x\in\partial D, \, t\in(0,\,T).
\end{equation*}
Then for any fixed $t\in(0,\,T)$ we have that $g_2(\cdot,\,t)\in H^{1/2}_\diamond(\partial D)$. Define $\varphi_1$ and $\varphi_2$ as the solutions to
\begin{equation}\label{bie_01}
\int_{\partial D} \frac{\varphi_1(y,\,t-c_0^{-1}|x-y|)}{4\pi |x-y|}\,ds(y) = g_1(t), \quad (x,\,t)\in (\partial D)_T,
\end{equation}
and
\begin{equation}\label{bie_02}
\int_{\partial D} \frac{\varphi_2(y,\,t-c_0^{-1}|x-y|)}{4\pi |x-y|}\,ds(y) = g_2(x,\,t), \quad (x,\,t)\in (\partial D)_T,
\end{equation}
respectively. Due to the linearity of the equation \eqref{bie}, we see that $\varphi = \varphi_1 + \varphi_2$.

Now, using the embedding $H^1(0,\,T)\hookrightarrow C[0,\,T]$, we have
\begin{eqnarray}\label{est_sigma1}
\|\varphi_2(\cdot,\,t)\|_{H^{-1/2}(\partial D)}&\lesssim& \| \varphi_2 \|_{H^1_{0,\sigma}(0,\,T;\,H^{-1/2}(\partial D))} \nonumber\\
&\lesssim &\left \| (\mathcal S_{\partial D})^{-1} \right\|_{\mathcal L \left(H^3_{0,\sigma}(0,\,T;\,H^{1/2}_\diamond(\partial D)),\,H^1_{0,\sigma}(0,\,T;\,H^{-1/2}(\partial D))\right)}\,\| g_2 \|_{H_{0,\sigma}^3(0,\,T;\,H^{1/2}(\partial D))} \nonumber\\
&\lesssim & \varepsilon\, \left \| (\mathcal S^\varepsilon_{\partial B})^{-1} \right\|_{\mathcal L \left(H^3_{0,\varepsilon\sigma}(0,\,T_\varepsilon;\,H^{1/2}_\diamond(\partial B)),\,H^1_{0,\varepsilon\sigma}(0,\,T_\varepsilon;\,H^{-1/2}(\partial B))\right)}\, \varepsilon^{1/2} \nonumber\\
&\lesssim & \varepsilon^{-1/2},
\end{eqnarray}
where we have used the fact that $\left \| (\mathcal S^\varepsilon_{\partial B})^{-1} \right\|_{\mathcal L \left(H^3_{0,\varepsilon\sigma}(0,\,T_\varepsilon;\,H^{1/2}_\diamond(\partial B)),\,H^1_{0,\varepsilon\sigma}(0,\,T_\varepsilon;\,H^{-1/2}(\partial B))\right)}$ is estimated by $O(1/\varepsilon^2)$ due to Lemma \ref{inv_SB}.

By the same derivation, we obtain
\begin{equation}\label{est_sigma2}
\|\varphi_1(\cdot,\,t)\|_{H^{-1/2}(\partial D)} \lesssim \| \varphi_1 \|_{H^1_{0,\sigma}(0,\,T;\,H^{-1/2}(\partial D))}  \lesssim \varepsilon^{-1/2}.
\end{equation}
Therefore, we have
\begin{equation}\label{est_sigma_3}
\|\varphi(\cdot,\,t)\|_{H^{-1/2}(\partial D)} \lesssim \| \varphi \|_{H^1_{0,\sigma}(0,\,T;\,H^{-1/2}(\partial D))}  \lesssim \varepsilon^{-1/2}.
\end{equation}
Similarly, we can also prove
\begin{equation}\label{est_sigmat_3}
\|\partial_t\varphi(\cdot,\,t)\|_{H^{-1/2}(\partial D)} \lesssim \| \partial_t\varphi \|_{H^1_{0,\sigma}(0,\,T;\,H^{-1/2}(\partial D))}  \lesssim \varepsilon^{-1/2}.
\end{equation}

\bigskip
We are now in a position to show the asymptotic behavior of the solution to \eqref{ibvp} with $M=1$.
\begin{theorem}\label{thS_main}
For $x\in\mathbb R^3\setminus\overline D$, with $x$ away from $D$, and $t\in (0,\,T)$, we have the following expansion:
\begin{equation}\label{S_main}
u^s(x,\,t) = -C_0\,\frac{\lambda(t-c_0^{-1}|x-z|-c_0^{-1}|z-z^*|)}{16\pi^2 |x-z|\,|z-z^*|} + O(\varepsilon^2) \quad \mathrm{as}\;\varepsilon \to 0,
\end{equation}
with the constant $C_0$ defined by
\begin{equation}\label{cap}
C_0:= \int_{\partial D} \varphi_0(y)\,ds(y),
\end{equation}
where $\varphi_0(x)$ is the solution to
\begin{equation}\label{bie_e}
\int_{\partial D} \frac{\varphi_0(y)}{4\pi |x-y|}\,ds(y) = 1, \quad x\in\partial D.
\end{equation}
\end{theorem}

{\bf Proof.} First, we rewrite the equation \eqref{bie} as
\begin{eqnarray*}
&&\int_{\partial D} \frac{\varphi(y,\,t)}{4\pi |x-y|}\,ds(y) + \int_{\partial D} \frac{\varphi(y,\,t-c_0^{-1}|x-y|) - \varphi(y,\,t)}{4\pi |x-y|}\,ds(y)\\
&=& \frac{-\lambda(t-c_0^{-1}|z-z^*|)}{4\pi |z-z^*|} + \left[ \frac{\lambda(t-c_0^{-1}|z-z^*|)}{4\pi |z-z^*|} - \frac{\lambda(t-c_0^{-1}|x-z^*|)}{4\pi |x-z^*|} \right], \quad (x,\,t)\in (\partial D)_T,\;z^*\not\in\overline D.
\end{eqnarray*}
Note that
\begin{eqnarray*}
\left|\int_{\partial D} \frac{\varphi(y,\,t-c_0^{-1}|x-y|) - \varphi(y,\,t)}{4\pi |x-y|}\,ds(y) \right|
&=& \left| \int_{\partial D} \frac{\partial_t \varphi(y,\,t^*)}{4\pi c_0} \,ds(y) \right| \\
&= & O(1)\, \langle \partial_t \varphi(\cdot,\,t^*),\,1\rangle_{\partial D} \\
&\lesssim & \|\partial_t\varphi(\cdot,\,t^*)\|_{H^{-1/2}(\partial D)}\,\|1\|_{H^{1/2}(\partial D)} \\
&\lesssim & O(\varepsilon),
\end{eqnarray*}
and
\begin{equation*}
\frac{\lambda(t-c_0^{-1}|z-z^*|)}{4\pi |z-z^*|} - \frac{\lambda(t-c_0^{-1}|x-z^*|)}{4\pi |x-z^*|} = O(\varepsilon), \quad (x,\,t)\in (\partial D)_T.
\end{equation*}
Then we conclude that
\begin{equation}\label{bie_1}
\int_{\partial D} \frac{\varphi(y,\,t)}{4\pi |x-y|}\,ds(y) = \frac{-\lambda(t-c_0^{-1}|z-z^*|)}{4\pi |z-z^*|} + O(\varepsilon),\quad (x,\,t)\in (\partial D)_T.
\end{equation}

Consider the equation
\begin{equation}\label{bie_2}
\int_{\partial D} \frac{\overline\varphi(y,\,t)}{4\pi |x-y|}\,ds(y) = \frac{-\lambda(t-c_0^{-1}|z-z^*|)}{4\pi |z-z^*|},\quad (x,\,t)\in (\partial D)_T.
\end{equation}
Then we have
\begin{equation*}\label{bie_4}
\overline\varphi(x,\,t) = \frac{-\lambda(t-c_0^{-1}|z-z^*|)}{4\pi |z-z^*|}\, \varphi_0(x),
\end{equation*}
and hence
\begin{equation}\label{bie_5}
\int_{\partial D}\overline\varphi(x,\,t)\,ds(x) = \frac{-\lambda(t-c_0^{-1}|z-z^*|)}{4\pi |z-z^*|}\,C_0.
\end{equation}
In addition, we obtain from \eqref{bie_1} and \eqref{bie_2} that
\begin{equation*}\label{bie_6}
\int_{\partial D} \frac{\varphi(y,\,t) - \overline \varphi(y,\,t)}{4\pi |x-y|}\,ds(y) = O(\varepsilon), \quad t\in(0,\,T),
\end{equation*}
and hence
\begin{equation}\label{bie_7}
\int_{\partial D} [\varphi(y,\,t) - \overline \varphi(y,\,t)]\,ds(y) = O(\varepsilon^2), \quad t\in(0,\,T).
\end{equation}

Now, for $(x,\,t)\in (\mathbb R^3\setminus\overline D)_T$, we derive
\begin{eqnarray*}\label{exp}
u^s(x,\,t) & = & \int_{\partial D} \frac{\varphi(y,\,t-c_0^{-1}|x-y|)}{4\pi |x-y|}\,ds(y) \nonumber\\
&=& \int_{\partial D} \frac{\varphi(y,\,t-c_0^{-1}|x-z|)}{4\pi |x-z|}\,ds(y) + \int_{\partial D}\left[ \frac{\varphi(y,\,t-c_0^{-1}|x-y|)}{4\pi |x-y|} - \frac{\varphi(y,\,t-c_0^{-1}|x-z|)}{4\pi |x-z|} \right]\,ds(y) \nonumber\\
&=& \int_{\partial D} \frac{\overline\varphi(y,\,t-c_0^{-1}|x-z|)}{4\pi |x-z|}\,ds(y) + \frac{1}{4\pi |x-z|} \int_{\partial D}[\varphi(y,\,t-c_0^{-1}|x-z|) - \overline\varphi(y,\,t-c_0^{-1}|x-z|)]\,ds(y) \nonumber \\
&& + \int_{\partial D}\varphi(y,\,t-c_0^{-1}|x-y|)\left[\frac{1}{4\pi |x-y|} - \frac{1}{4\pi |x-z|} \right]\,ds(y) \nonumber \\
&& + \int_{\partial D} \frac{1}{4\pi |x-z|}\left[\varphi(y,\,t-c_0^{-1}|x-y|) - \varphi(y,\,t-c_0^{-1}|x-z|) \right]\,ds(y) \nonumber\\
&=& -C_0\,\frac{\lambda(t-c_0^{-1}|x-z|-c_0^{-1}|z-z^*|)}{16\pi^2 |x-z|\,|z-z^*|} + O(\varepsilon^2) + O(\varepsilon) \int_{\partial D}\varphi(y,\,t-c_0^{-1}|x-y|)\,ds(y) \nonumber\\
&& + O(\varepsilon) \int_{\partial D} \partial_t \varphi(y,\,t^*)\,ds(y) \quad (\mbox{as } x \mbox{ is away from } D)  \nonumber\\
&=& -C_0\,\frac{\lambda(t-c_0^{-1}|x-z|-c_0^{-1}|z-z^*|)}{16\pi^2 |x-z|\,|z-z^*|} + O(\varepsilon^2).
\end{eqnarray*}
The proof is complete. \hfill $\Box$

\section{Proof of Theorem \ref{Main}: the multiple holes case}\label{multiple}
\setcounter{equation}{0}

In this section, we give a rigorous justification of the asymptotic expansion for the solution to \eqref{ibvp} as $\varepsilon \ll 1$ and prove 
the unique solvability of the linear algebraic system \eqref{M1Mq_i} for the case of multiple holes.

We express the solution to \eqref{ibvp} as a retarded single-layer potential
\begin{equation}\label{Mlayer}
u^s(x,\,t)=\sum_{j=1}^M\int_{\partial D_j} \frac{\varphi_j(y,\,t-c_0^{-1}|x-y|)}{4\pi |x-y|}\,ds(y), \quad (x,\,t)\in(\mathbb R^3\setminus\overline D)_T,
\end{equation}
where $\varphi_j$'s are causal densities to be determined. Then we have
\begin{eqnarray}\label{Mbie}
&&\int_{\partial D_i} \frac{\varphi_i(y,\,t-c_0^{-1}|x-y|)}{4\pi |x-y|}\,ds(y) + \sum_{\substack{j=1\\j\neq i}}^M \int_{\partial D_j} \frac{\varphi_j(y,\,t-c_0^{-1}|x-y|)}{4\pi |x-y|}\,ds(y) \nonumber\\
& = & \frac{-\lambda(t-c_0^{-1}|x-z^*|)}{4\pi |x-z^*|}, \quad (x,\,t)\in (\partial D_i)_T,\;i=1,\,2,\,\cdots,\,M.
\end{eqnarray}
For convenience, we define
\begin{equation}\label{S_ij}
\mathcal S_{ji}[\varphi_j](x,\,t):=\int_{\partial D_j} \frac{\varphi_j(y,\,t-c_0^{-1}|x-y|)}{4\pi |x-y|}\,ds(y), \quad (x,\,t)\in (\partial D_i)_T,
\end{equation}
and rewrite \eqref{Mbie} as
\begin{equation}\label{Mbie2}
\varphi_i + \sum_{j\not= i}\mathcal S_{ii}^{-1} \mathcal S_{ji}[\varphi_j] = - \mathcal S_{ii}^{-1}[u^i] \quad \mathrm{on}\;(\partial D_i)_T.
\end{equation}
The unique solvability of \eqref{Mbie} can be shown in a standard way. To proceed, we first prove the following {\it a priori} estimate of the densities.

\begin{lemma}\label{thM_est_phi}
Define $d_{i*}=\mathrm{dist}(z^*,\,D_i)$ and $\mathcal C_n:= \max_{t\in [0,\,T]}|\lambda^{(n)}(t)|$. Assume that $\sum_{n=0}^{+\infty}\mathcal C_n$ is convergent and denote $\mathcal C:=\sum_{n=0}^{+\infty}\mathcal C_n$. Then, under the condition
\begin{equation}\label{M_cond}
\varepsilon^{2-r} \max_{1\leq i \leq M}\sum_{j\not=i} d_{ij}^{-2} < 1,
\end{equation} 
we have
\begin{eqnarray}\label{M_est_phi}
\sum_{n=0}^{+\infty} \left \| \frac{\partial^n\varphi_i(y,\,t)}{\partial t^n} \right\|_{H^r_{0,\sigma}(0,\,T;\,H^{-1/2}(\partial D_i))}
& \lesssim & \left( 1 -  \varepsilon^{2-r} \max_{1\leq i \leq M}\sum_{j\not=i} d_{ij}^{-2}\right)^{-1}\,\left(\mathcal C\,\varepsilon^{1/2-r}\right) \max_{1\leq i \leq M} d_{i*}^{-1} \nonumber\\
&=& O(\varepsilon^{1/2-r}), \qquad r=0,\,1,\;i=1,\,2,\,\cdots,\,M.
\end{eqnarray}
\end{lemma}

{\bf Proof.} First, by the same argument as for estimating $(\mathcal S_{\partial D})^{-1}$ in Section \ref{single}, we can derive the estimate
\begin{eqnarray}\label{MS_ii_inv}
& & \left \| \mathcal S_{ii}^{-1} \right\|_{\mathcal L \left(H^{r+2}_{0,\sigma}(0,\,T;\,H^{1/2}(\partial D_i)),\,H^r_{0,\sigma}(0,\,T;\,H^{-1/2}(\partial D_i))\right)} \nonumber \\
& \leq &\left \| (\mathcal S^\varepsilon_{\partial B_i})^{-1} \right\|_{\mathcal L \left(H^{r+2}_{0,\varepsilon\sigma}(0,\,T_\varepsilon;\,H^{1/2}(\partial B_i)),\,H^r_{0,\varepsilon\sigma}(0,\,T_\varepsilon;\,H^{-1/2}(\partial B_i))\right)} \lesssim \varepsilon^{-(r+1)}.
\end{eqnarray}
Next, we show the estimate of $\left \| \mathcal S_{ji}[\varphi_j] \right\|_{H^{r+2}_{0,\sigma}(0,\,T;\,H^{1/2}(\partial D_i))}$ for $j\not=i$.
Since
\begin{equation}\label{MS_ji_1}
\left \| \mathcal S_{ji}[\varphi_j] \right\|^2_{H^{r+2}_{0,\sigma}(0,\,T;\,H^{1/2}(\partial D_i))} \lesssim C(T) \int_0^T e^{-2\sigma t} \left\| \frac{\partial^{r+2} \mathcal S_{ji}[\varphi_j]}{\partial t^{r+2}} \right\|^2_{H^{1/2}(\partial D_i)}\,dt,
\end{equation}
it suffices to estimate $\left\| \frac{\partial^{r+2} \mathcal S_{ji}[\varphi_j]}{\partial t^{r+2}} \right\|^2_{H^{1/2}(\partial D_i)}$. By the definition  \eqref{defn_Hnorm1} of the $H^{1/2}$-norm, let us estimate $\left\| \frac{\partial^{r+2} \mathcal S_{ji}[\varphi_j]}{\partial t^{r+2}} \right\|^2_{H^1(D_i)}$. Note that for $x\in D_i$ we have
\begin{eqnarray*}
\left| \frac{\partial^{r+2} \mathcal S_{ji}[\varphi_j]}{\partial t^{r+2}} \right| & = & \left| \frac{\partial^{r+2}}{\partial t^{r+2}}\int_{\partial D_j} \frac{\varphi_j(y,\,t-c_0^{-1}|x-y|)}{4\pi |x-y|}\,ds(y) \right| \nonumber\\
&=& \left| \int_{\partial D_j} \frac{1}{4\pi |x-y|}\frac{\partial^{r+2}\varphi_j(y,\,t-c_0^{-1}|x-y|)}{\partial t^{r+2}}\,ds(y) \right| \nonumber\\
&=& \left| \int_{\partial D_j} \frac{1}{4\pi |x-y|}\left[\frac{\partial^{r+2}\varphi_j(y,\,t)}{\partial t^{r+2}} - \frac{\partial^{r+3}\varphi_j(y,\,t-\theta_jc_0^{-1}|x-y|)}{\partial t^{r+3}}\left(c_0^{-1}|x-y|\right) \right]\,ds(y) \right| \nonumber\\
&\lesssim & d_{ij}^{-1}\left|\left\langle \frac{\partial^{r+2}\varphi_j(y,\,t)}{\partial t^{r+2}},\,1\right\rangle\right| + \left|\left\langle \frac{\partial^{r+3}\varphi_j(y,\,t-\theta_jc_0^{-1}|x-y|)}{\partial t^{r+3}},\,1\right\rangle\right| \nonumber\\
&\lesssim & d_{ij}^{-1} \left\| \frac{\partial^{r+2}\varphi_j(y,\,t)}{\partial t^{r+2}} \right\|_{H^{-1/2}(\partial D_j)}\, \| 1 \|_{H^{1/2}(\partial D_j)} \nonumber \\
&& \qquad \qquad \qquad + \left\| \frac{\partial^{r+3}\varphi_j(y,\,t-\theta_jc_0^{-1}|x-y|)}{\partial t^{r+3}} \right\|_{H^{-1/2}(\partial D_j)}\, \left\| 1 \right\|_{H^{1/2}(\partial D_j)} \nonumber\\
&\lesssim & d_{ij}^{-1}\,\varepsilon^{3/2}\,\left\| \frac{\partial^{r+2}\varphi_j(y,\,t)}{\partial t^{r+2}} \right\|_{H^{-1/2}(\partial D_j)} + \varepsilon^{3/2}\,\left\| \frac{\partial^{r+3}\varphi_j(y,\,t-\theta_jc_0^{-1}|x-y|)}{\partial t^{r+3}} \right\|_{H^{-1/2}(\partial D_j)},
\end{eqnarray*}
where $\theta_j\in(0,\,1)$. Then it follows that
\begin{eqnarray*}\label{MS_ji_2}
&&\int_{D_i} \left| \frac{\partial^{r+2} \mathcal S_{ji}[\varphi_j]}{\partial t^{r+2}} \right|^2\,ds(x) \nonumber\\
&\lesssim & d_{ij}^{-2}\,\varepsilon^6\,\left\| \frac{\partial^{r+2}\varphi_j(y,\,t)}{\partial t^{r+2}} \right\|^2_{H^{-1/2}(\partial D_j)} + \varepsilon^3\,\int_{D_i}\left\| \frac{\partial^{r+3}\varphi_j(y,\,t-\theta_jc_0^{-1}|x-y|)}{\partial t^{r+3}} \right\|^2_{H^{-1/2}(\partial D_j)}\,ds(x).
\end{eqnarray*}
Similarly, we can also estimate $\left\|\nabla_x \frac{\partial^{r+2} \mathcal S_{ji}[\varphi_j]}{\partial t^{r+2}} \right\|^2_{L^2(D_i)}$ as
\begin{eqnarray*}\label{MS_ji_2_g}
&&\int_{D_i} \left|\nabla_x \frac{\partial^{r+2} \mathcal S_{ji}[\varphi_j]}{\partial t^{r+2}} \right|^2\,ds(x) \nonumber\\
&\lesssim & d_{ij}^{-4}\,\varepsilon^6\,\left\| \frac{\partial^{r+2}\varphi_j(y,\,t)}{\partial t^{r+2}} \right\|^2_{H^{-1/2}(\partial D_j)} + d_{ij}^{-2}\, \varepsilon^3\,\int_{D_i}\left\| \frac{\partial^{r+3}\varphi_j(y,\,t-\theta_jc_0^{-1}|x-y|)}{\partial t^{r+3}} \right\|^2_{H^{-1/2}(\partial D_j)}\,ds(x) \nonumber\\
&& + d_{ij}^{-2}\,\varepsilon^6\,\left\| \frac{\partial^{r+3}\varphi_j(y,\,t)}{\partial t^{r+3}} \right\|^2_{H^{-1/2}(\partial D_j)} + \varepsilon^3\,\int_{D_i}\left\| \frac{\partial^{r+4}\varphi_j(y,\,t-\tilde\theta_jc_0^{-1}|x-y|)}{\partial t^{r+4}} \right\|^2_{H^{-1/2}(\partial D_j)}\,ds(x)
\end{eqnarray*}
with $\theta_j,\,\tilde\theta_j\in(0,\,1)$. So we finally have
\begin{eqnarray*}\label{MS_ji_22}
&&\left\|\frac{\partial^{r+2} \mathcal S_{ji}[\varphi_j]}{\partial t^{r+2}} \right\|^2_{H^{1/2}(\partial D_i)}\leq \left\|\frac{\partial^{r+2} \mathcal S_{ji}[\varphi_j]}{\partial t^{r+2}} \right\|^2_{H^1(D_i)}\nonumber\\
&\lesssim &d_{ij}^{-4}\,\varepsilon^6\,\left\| \frac{\partial^{r+2}\varphi_j(y,\,t)}{\partial t^{r+2}} \right\|^2_{H^{-1/2}(\partial D_j)} + d_{ij}^{-2}\varepsilon^3\,\int_{D_i}\left\| \frac{\partial^{r+3}\varphi_j(y,\,t-\theta_jc_0^{-1}|x-y|)}{\partial t^{r+3}} \right\|^2_{H^{-1/2}(\partial D_j)}\,ds(x) \nonumber\\
& & + d_{ij}^{-2}\,\varepsilon^6\,\left\| \frac{\partial^{r+3}\varphi_j(y,\,t)}{\partial t^{r+3}} \right\|^2_{H^{-1/2}(\partial D_j)} + \varepsilon^3\,\int_{D_i}\left\| \frac{\partial^{r+4}\varphi_j(y,\,t-\tilde\theta_jc_0^{-1}|x-y|)}{\partial t^{r+4}} \right\|^2_{H^{-1/2}(\partial D_j)}\,ds(x).
\end{eqnarray*}
Thus, we derive
\begin{eqnarray}\label{MS_ji_3}
& & \left \| \mathcal S_{ji}[\varphi_j] \right\|^2_{H^{r+2}_{0,\sigma}(0,\,T;\,H^{1/2}(\partial D_i))} \nonumber\\
&\lesssim &  d_{ij}^{-4}\,\varepsilon^6\,\int_0^T e^{-2\sigma t} \,\left\| \frac{\partial^{r+2}\varphi_j(y,\,t)}{\partial t^{r+2}} \right\|^2_{H^{-1/2}(\partial D_j)} \,dt + d_{ij}^{-2}\,\varepsilon^6\,\int_0^T e^{-2\sigma t} \,\left\| \frac{\partial^{r+3}\varphi_j(y,\,t)}{\partial t^{r+3}} \right\|^2_{H^{-1/2}(\partial D_j)} \,dt \nonumber \\
&& \qquad + d_{ij}^{-2}\,\varepsilon^3  \int_0^T e^{-2\sigma t} \,\int_{D_i}\left\| \frac{\partial^{r+3}\varphi_j(y,\,t-\theta_jc_0^{-1}|x-y|)}{\partial t^{r+3}} \right\|^2_{H^{-1/2}(\partial D_j)}\,ds(x) \,dt \nonumber\\
&& \qquad + \varepsilon^3  \int_0^T e^{-2\sigma t} \,\int_{D_i}\left\| \frac{\partial^{r+4}\varphi_j(y,\,t-\tilde\theta_jc_0^{-1}|x-y|)}{\partial t^{r+4}} \right\|^2_{H^{-1/2}(\partial D_j)}\,ds(x) \,dt \nonumber\\
&\lesssim & d_{ij}^{-4}\,\varepsilon^6 \,\left\| \frac{\partial^2\varphi_j}{\partial t^2} \right\|^2_{H^r_{0,\sigma}(0,\,T;\,H^{-1/2}(\partial D_j))} + d_{ij}^{-2}\,\varepsilon^6 \,\left\| \frac{\partial^3\varphi_j}{\partial t^3} \right\|^2_{H^r_{0,\sigma}(0,\,T;\,H^{-1/2}(\partial D_j))} \nonumber\\
&& \qquad + d_{ij}^{-2}\,\varepsilon^3 \int_{D_i} \int_{-\theta_jc_0^{-1}|x-y|}^{T-\theta_jc_0^{-1}|x-y|} e^{-2\sigma (\eta+\theta_jc_0^{-1}|x-y|)} \,\left\| \frac{\partial^{r+3}\varphi_j(y,\,\eta)}{\partial \eta^{r+3}} \right\|^2_{H^{-1/2}(\partial D_j)}\,d\eta ds(x)\nonumber\\
&& \qquad + \varepsilon^3 \int_{D_i} \int_{-\tilde\theta_jc_0^{-1}|x-y|}^{T-\tilde\theta_jc_0^{-1}|x-y|} e^{-2\sigma (\eta+\tilde\theta_jc_0^{-1}|x-y|)} \,\left\| \frac{\partial^{r+4}\varphi_j(y,\,\eta)}{\partial \eta^{r+4}} \right\|^2_{H^{-1/2}(\partial D_j)}\,d\eta ds(x)\nonumber\\
&\lesssim & d_{ij}^{-4}\,\varepsilon^6 \,\left\| \frac{\partial^2\varphi_j}{\partial t^2} \right\|^2_{H^r_{0,\sigma}(0,\,T;\,H^{-1/2}(\partial D_j))} + d_{ij}^{-2}\,\varepsilon^6 \,\left\| \frac{\partial^3\varphi_j}{\partial t^3} \right\|^2_{H^r_{0,\sigma}(0,\,T;\,H^{-1/2}(\partial D_j))} \nonumber\\
&& \qquad + d_{ij}^{-2}\,\varepsilon^3 \int_{D_i} \int_0^T e^{-2\sigma \eta} \,\left\| \frac{\partial^{r+3}\varphi_j(y,\,\eta)}{\partial \eta^{r+3}} \right\|^2_{H^{-1/2}(\partial D_j)}\,d\eta ds(x)\nonumber\\
&& \qquad  + \varepsilon^3 \int_{D_i} \int_0^T e^{-2\sigma \eta} \,\left\| \frac{\partial^{r+4}\varphi_j(y,\,\eta)}{\partial \eta^{r+4}} \right\|^2_{H^{-1/2}(\partial D_j)}\,d\eta ds(x).
\end{eqnarray}
It implies that
\begin{eqnarray}\label{MS_ji_4}
\left \| \mathcal S_{ji}[\varphi_j] \right\|_{H^{r+2}_{0,\sigma}(0,\,T;\,H^{1/2}(\partial D_i))} 
&\lesssim & \varepsilon^3 d_{ij}^{-2}\,\left\| \frac{\partial^2\varphi_j}{\partial t^2} \right\|_{H^r_{0,\sigma}(0,\,T;\,H^{-1/2}(\partial D_j))} + \varepsilon^3\,d_{ij}^{-1} \left\| \frac{\partial^3\varphi_j}{\partial t^3} \right\|_{H^r_{0,\sigma}(0,\,T;\,H^{-1/2}(\partial D_j))}\nonumber\\
&& + \varepsilon^3\,\left\| \frac{\partial^4\varphi_j}{\partial t^4} \right\|_{H^r_{0,\sigma}(0,\,T;\,H^{-1/2}(\partial D_j))}.
\end{eqnarray}
Then, from \eqref{Mbie2} and \eqref{MS_ii_inv}, we derive
\begin{eqnarray}\label{MS_ji_5}
&&\left \| \varphi_i \right\|_{H^r_{0,\sigma}(0,\,T;\,H^{-1/2}(\partial D_i))}\nonumber\\
&\lesssim & \varepsilon^{2-r} \sum_{j\not=i} d_{ij}^{-2}\,\left\| \frac{\partial^2\varphi_j}{\partial t^2} \right\|_{H^r_{0,\sigma}(0,\,T;\,H^{-1/2}(\partial D_j))} + \varepsilon^{2-r} \sum_{j\not=i} d_{ij}^{-1}  \left\| \frac{\partial^3\varphi_j}{\partial t^3} \right\|_{H^r_{0,\sigma}(0,\,T;\,H^{-1/2}(\partial D_j))} \nonumber\\
&& + \varepsilon^{2-r} \sum_{j\not=i} \left\| \frac{\partial^4\varphi_j}{\partial t^4} \right\|_{H^r_{0,\sigma}(0,\,T;\,H^{-1/2}(\partial D_j))} + \varepsilon^{1/2-r}\, \mathcal C_{r+2}\, d_{i*}^{-1}.
\end{eqnarray}

Using the same argument, we can estimate $\left \| \frac{\partial^n\varphi_i}{\partial t^n} \right\|_{H^r_{0,\sigma}(0,\,T;\,H^{-1/2}(\partial D_i))} $ for $n=2,\,3,\,\cdots$ as
\begin{eqnarray}\label{MS_ji_6}
&&\left \|  \frac{\partial^n\varphi_i}{\partial t^n} \right\|_{H^r_{0,\sigma}(0,\,T;\,H^{-1/2}(\partial D_i))} \nonumber\\
&\lesssim & \varepsilon^{2-r} \sum_{j\not=i} d_{ij}^{-2}\,\left\| \frac{\partial^{n+2}\varphi_j}{\partial t^{n+2}} \right\|_{H^r_{0,\sigma}(0,\,T;\,H^{-1/2}(\partial D_j))} + \varepsilon^{2-r} \sum_{j\not=i} d_{ij}^{-1}  \left\| \frac{\partial^{n+3}\varphi_j}{\partial t^{n+3}} \right\|_{H^r_{0,\sigma}(0,\,T;\,H^{-1/2}(\partial D_j))} \nonumber\\
&&\qquad \qquad + \varepsilon^{2-r} \sum_{j\not=i} \left\| \frac{\partial^{n+4}\varphi_j}{\partial t^{n+4}} \right\|_{H^r_{0,\sigma}(0,\,T;\,H^{-1/2}(\partial D_j))} + \varepsilon^{1/2-r}\,\mathcal C_{r+2+n}\, d_{i*}^{-1}.
\end{eqnarray}
Let
\begin{equation*}
\mathcal A_i:= \sum_{n=0}^{+\infty} \left \| \frac{\partial^n\varphi_i}{\partial t^n} \right\|_{H^r_{0,\sigma}(0,\,T;\,H^{-1/2}(\partial D_i))}.
\end{equation*}
Then we have
\begin{eqnarray*}\label{MS_ji_7}
\mathcal A_i &\lesssim & 3\, \varepsilon^{2-r} \sum_{j\not=i} d_{ij}^{-2}\, \mathcal A_j + \varepsilon^{1/2-r}d_{i*}^{-1}\sum_{n=0}^{+\infty}\mathcal C_{r+2+n} \nonumber\\
&\lesssim& \varepsilon^{2-r}\,\left(\max_{1\leq j \leq M}\mathcal A_j\right)\, \sum_{j\not=i} d_{ij}^{-2} + \mathcal C\,\varepsilon^{1/2-r}\,\max_{1\leq i \leq M} d_{i*}^{-1}.
\end{eqnarray*}
Under the condition \eqref{M_cond}, we get
\begin{eqnarray*}\label{MS_ji_9}
\max_{1\leq j \leq M}\,\mathcal A_j \lesssim  \left( 1 -  \varepsilon^{2-r} \max_{1\leq i \leq M}\sum_{j\not=i} d_{ij}^{-2}\right)^{-1}\,\left(\mathcal C\,\varepsilon^{1/2-r}\right) \max_{1\leq i \leq M} d_{i*}^{-1}.
\end{eqnarray*}
The proof is complete. \hfill $\Box$

\begin{remark}\label{reM_est}
Using the embedding $H^1(0,\,T)\hookrightarrow C[0,\,T]$ and taking $r=1$ in \eqref{M_est_phi}, we have the pointwise estimate
\begin{equation}\label{M_est_phi_p}
\left| \frac{\partial^n \varphi_i(\cdot,\,t)}{\partial t^n} \right|_{H^{-1/2}(\partial D_i)} = O\left(\varepsilon^{-1/2}\right),\quad t\in(0,\,T),\; i=1,\,\cdots,\,M,\; n=0,\,1,\,\cdots.
\end{equation}
\end{remark}

\bigskip
Set
\begin{equation}\label{Mv_i}
v_i(t) := \int_{\partial D_i} \varphi_i(y,\,t)\,ds(y), \quad t\in(0,\,T),\; i=1,\,2,\,\cdots,\,M.
\end{equation}
Then \eqref{Mbie} can be rewritten as
\begin{eqnarray}\label{Mbie1}
&&\int_{\partial D_i} \frac{\varphi_i(y,\,t)}{4\pi |x-y|}\,ds(y) +  \sum_{\substack{j=1\\j\neq i}}^M \frac{v_j(t-c_0^{-1}|z_i-z_j|)}{4\pi |z_i-z_j|} \nonumber\\
&=& \frac{-\lambda(t-c_0^{-1}|z_i-z^*|)}{4\pi |z_i-z^*|} + E_i^{(1)} + E_i^{(2)} + E_i^{(3)}, \quad (x,\,t)\in (\partial D_i)_T,
\end{eqnarray}
where
\begin{eqnarray}\label{ME_1}
E_i^{(1)} &:=& \int_{\partial D_i} \frac{\varphi_i(y,\,t) - \varphi_i(y,\,t-c_0^{-1}|x-y|)}{4\pi |x-y|}\,ds(y) \nonumber\\
&& +  \sum_{\substack{j=1\\j\neq i}}^M \int_{\partial D_j} \frac{\varphi_j(y,\,t-c_0^{-1}|z_i-z_j|)-\varphi_j(y,\,t-c_0^{-1}|x-y|)}{4\pi |x-y|}\,ds(y),
\end{eqnarray}
\begin{equation}\label{ME_2}
E_i^{(2)} := \sum_{\substack{j=1\\j\neq i}}^M  \int_{\partial D_j} \left( \frac{1}{4\pi |z_i-z_j|} - \frac{1}{4\pi |x-y|}\right)\,\varphi_j(y,\,t-c_0^{-1}|z_i-z_j|)\,ds(y),
\end{equation}
and
\begin{equation}\label{ME_3}
E_i^{(3)} := \frac{\lambda(t-c_0^{-1}|z_i-z^*|)}{4\pi |z_i-z^*|} - \frac{\lambda(t-c_0^{-1}|x-z^*|)}{4\pi |x-z^*|}.
\end{equation}

\bigskip
In the following, let us estimate $E_i^{(1)},\,E_i^{(2)}$ and $E_i^{(3)}$. Since
\begin{eqnarray*}\label{ME_1_1}
&&\left|\int_{\partial D_i} \frac{\varphi_i(y,\,t) - \varphi_i(y,\,t-c_0^{-1}|x-y|)}{4\pi |x-y|}\,ds(y) \right| \nonumber\\
&=& \frac{1}{4\pi} \left| \int_{\partial D_i} \frac{\partial  \varphi_i(y,\,t^*)}{\partial t}\,ds(y) \right| \nonumber\\
&\leq &\frac{1}{4\pi} \int_{\partial D_i} \max_{t\in[0,\,T]} \left|\frac{\partial  \varphi_i(y,\,t)}{\partial t}\right|\,ds(y) \nonumber\\
&\lesssim & \left\langle\max_{t\in[0,\,T]} \left|\frac{\partial  \varphi_i(y,\,t)}{\partial t}\right|,\,1\right\rangle \lesssim  \varepsilon.
\end{eqnarray*}
For $i\not=j$, we have
\begin{eqnarray*}\label{ME_1_3}
& &\left|\int_{\partial D_j} \frac{\varphi_j(y,\,t-c_0^{-1}|z_i-z_j|)-\varphi_j(y,\,t-c_0^{-1}|x-y|)}{4\pi |x-y|}\,ds(y)\right| \\
&=& \left| \int_{\partial D_j} \frac{\varphi_j(y,\,t-c_0^{-1}|z_i-z_j|)-\varphi_j(y,\,t-c_0^{-1}|x-z_j|)}{4\pi |x-y|}\,ds(y)\right|\\
& & + \left|\int_{\partial D_j} \frac{\varphi_j(y,\,t-c_0^{-1}|x-z_j|)-\varphi_j(y,\,t-c_0^{-1}|y-z_j|)}{4\pi |x-y|}\,ds(y)\right|\\
&=&
O\left(\varepsilon^2\,d_{ij}^{-1} + \varepsilon \right).
\end{eqnarray*}
So we obtain
\begin{equation}\label{ME_1_4}
E_i^{(1)} = O(\varepsilon) + O\big(\varepsilon^2 \sum_{j\not=i} d_{ij}^{-1}\big).
\end{equation}
Due to the estimate
\begin{equation*}
\left|\frac{1}{|z_i-z_j|} - \frac{1}{|x-y|}\right|=\frac{\big||x-y|-|z_i-z_j|\big|}{|z_i-z_j|\,|x-y|}\leq 2a\,d_{ij}^{-2}, \quad x\in\partial D_i,\,y\in\partial D_j,\,i\not=j,
\end{equation*}
we can easily see
\begin{equation}\label{ME_2_1}
E_i^{(2)} = O\big(\varepsilon^2 \sum_{j\not=i} d_{ij}^{-2}\big).
\end{equation}
In addition, it can be easily deduced that
\begin{equation}\label{ME_31}
E_i^{(3)} = O(\varepsilon).
\end{equation}

Let $\overline\varphi_i(x,\,t)$ be the solution to
\begin{equation}\label{Mbphi}
\int_{\partial D_i} \frac{\overline \varphi_i(y,\,t)}{4\pi |x-y|}\,ds(y) = -\sum_{\substack{j=1\\j\neq i}}^M \frac{v_j(t-c_0^{-1}|z_i-z_j|)}{4\pi |z_i - z_j|} - \frac{\lambda(t-c_0^{-1}|z_i - z^*|)}{4\pi |z_i - z^*|}, \quad x\in\partial D_i,\;t\in(0,\,T).
\end{equation}
Define
\begin{equation}\label{Mbv_i0}
\overline v_i(t):= \int_{\partial D_i} \overline \varphi_i(y,\,t)\,ds(y).
\end{equation}
Then we have
\begin{equation}\label{Mbv_i1}
\overline v_i(t) = - C_i \sum_{\substack{j=1\\j\neq i}}^M \frac{v_j(t-c_0^{-1}|z_i-z_j|)}{4\pi |z_i - z_j|} - C_i\frac{\lambda(t-c_0^{-1}|z_i - z^*|)}{4\pi |z_i - z^*|}, \quad t\in(0,\,T).
\end{equation}

From \eqref{Mbie1} and \eqref{Mbphi}, we obtain for $(x,\,t)\in (\partial D_i)_T$ that
\begin{equation*}
\int_{\partial D_i} \frac{\varphi_i(y,\,t) - \overline \varphi_i(y,\,t)}{4\pi |x-y|}\,ds(y)
= E_i^{(1)} + E_i^{(2)} + E_i^{(3)}= O(\varepsilon) + O\big(\varepsilon^2 \sum_{j\not=i} d_{ij}^{-2}\big)=:E_i^{(4)},
\end{equation*}
and hence
\begin{equation}\label{Mbphi1}
\int_{\partial D_i} \left[\varphi_i(y,\,t) - \overline \varphi_i(y,\,t)\right]\,ds(y) = C_i\, E_i^{(4)}.
\end{equation}
By \eqref{Mbv_i1} and \eqref{Mbphi1}, we get
\begin{equation}\label{Mv_i_1}
C_i^{-1} v_i(t) + \sum_{\substack{j=1\\j\not= i}}^M \frac{v_j(t-c_0^{-1}|z_i-z_j|)}{4\pi |z_i-z_j|} = - \frac{\lambda(t-c_0^{-1}|z_i-z^*|)}{4\pi |z_i-z^*|} + E_i^{(4)}, \quad t\in (0,\,T),\;i=1,\,2,\,\cdots,\,M.
\end{equation}

To proceed, we show the invertibility of this linear algebraic system.

\begin{lemma}\label{Mlinear}
If
\begin{equation}\label{Mcond}
C \max_{1\leq i \leq M} \sum_{j\not= i} \frac{1}{4\pi |z_i-z_j|} < 1
\end{equation}
with $C:=\max_{1\leq j \leq M} C_j$, then the linear algebraic system
\begin{equation}\label{Mlinear_1}
q_i(t) + \sum_{\substack{j=1\\j\not= i}}^M \frac{C_jq_j(t-c_0^{-1}|z_i-z_j|)}{4\pi |z_i-z_j|} = f_i(t), \quad t\in(0,\,T),\;i=1,\,2,\,\cdots,\,M
\end{equation}
is uniquely solvable in $H_0^1(0,\,T)$. Moreover, we have the estimate
\begin{equation}\label{Mlinear_2}
\left( \sum_{i=1}^M |q_i(t)|^2 \right)^{1/2} \leq \left( 1 - C \max_{1\leq i \leq M} \sum_{j\not=i}\frac{1}{4\pi |z_i-z_j|}\right)^{-1} \left( \sum_{i=1}^M \Vert f_i\Vert_{H^1(0,\,T)}^2 \right)^{1/2}.
\end{equation}
\end{lemma}

{\bf Proof.} Multiplying \eqref{Mlinear_1} by $q_i(t)$ and taking integration with respect to $t$, we have
\begin{equation*}
\sum_{i=1}^M \int_0^T q^2_i(t)\,dt + \sum_{i=1}^M \sum_{\substack{j=1\\j\not= i}}^M \frac{C_j\int_0^T q_i(t)q_j(t-c_0^{-1}|z_i-z_j|)\,dt}{4\pi |z_i-z_j|} = \sum_{i=1}^M \int_0^T f_i(t)\,q_i(t)\,dt.
\end{equation*}
Since
\begin{eqnarray*}
\sum_{i=1}^M \sum_{\substack{j=1\\j\not= i}}^M \frac{1}{4\pi |z_i-z_j|}\int_0^T q_i(t)q_j(t)\,dt &\leq& \sum_{i=1}^M \sum_{\substack{j=1\\j\not= i}}^M \frac{1}{4\pi |z_i-z_j|}\left(\int_0^T q_i^2(t)\,dt\right)^{1/2}\,\left(\int_0^T q_j^2(t-c_0^{-1}|z_i-z_j|)\,dt\right)^{1/2}\\ 
&\leq& \sum_{i=1}^M \sum_{\substack{j=1\\j\not= i}}^M \frac{\Vert q_i\Vert_{L^2(0,\,T)}\,\Vert q_j\Vert_{L^2(0,\,T)}}{4\pi |z_i-z_j|}\\
&\leq& \sum_{i=1}^M \sum_{\substack{j=1\\j\not= i}}^M \frac{\Vert q_i\Vert_{L^2(0,\,T)}^2}{4\pi |z_i-z_j|},
\end{eqnarray*}
it follows that
\begin{equation*}
\sum_{i=1}^M \Vert q_i\Vert_{L^2(0,\,T)}^2 - C \sum_{i=1}^M \sum_{\substack{j=1\\j\not= i}}^M \frac{\Vert q_i\Vert_{L^2(0,\,T)}^2}{4\pi |z_i-z_j|} \leq  \left(\sum_{i=1}^M \Vert f_i\Vert_{L^2(0,\,T)}^2\right)^{1/2}\,\left(\sum_{i=1}^M \Vert q_i\Vert_{L^2(0,\,T)}^2\right)^{1/2}.
\end{equation*}
So we get
\begin{equation*}\label{Mlinear_3}
\left( 1 - C \max_{1\leq i \leq M} \sum_{j\not=i} \frac{1}{4\pi |z_i-z_j|}\right)\, \sum_{i=1}^M \Vert q_i\Vert_{L^2(0,\,T)}^2  \leq  \left(\sum_{i=1}^M \Vert f_i\Vert_{L^2(0,\,T)}^2\right)^{1/2}\,\left(\sum_{i=1}^M \Vert q_i\Vert_{L^2(0,\,T)}^2\right)^{1/2}.
\end{equation*}
Due to the condition \eqref{Mcond}, we have
\begin{equation*}
\left(\sum_{i=1}^M \Vert q_i\Vert_{L^2(0,\,T)}^2\right)^{1/2}  \leq \left( 1 - C \max_{1\leq i \leq M} \sum_{j\not=i} \frac{1}{4\pi |z_i-z_j|}\right)^{-1} \left(\sum_{i=1}^M \Vert f_i\Vert_{L^2(0,\,T)}^2\right)^{1/2}.
\end{equation*}
Take the derivative with respect to $t$ for \eqref{Mlinear_1} and use the same argument as above for $q^\prime_i(t)$. Then we have
\begin{equation*}
\left(\sum_{i=1}^M \Vert q^\prime_i\Vert_{L^2(0,\,T)}^2\right)^{1/2}  \leq \left( 1 - C \max_{1\leq i \leq M} \sum_{j\not=i} \frac{1}{4\pi |z_i-z_j|}\right)^{-1} \left(\sum_{i=1}^M \Vert f^\prime_i\Vert_{L^2(0,\,T)}^2\right)^{1/2}.
\end{equation*}
So we obtain
\begin{equation*}
\left(\sum_{i=1}^M \Vert q_i\Vert_{H^1(0,\,T)}^2\right)^{1/2}  \leq \left( 1 - C \max_{1\leq i \leq M} \sum_{j\not=i} \frac{1}{4\pi |z_i-z_j|}\right)^{-1} \left(\sum_{i=1}^M \Vert f_i\Vert_{H^1(0,\,T)}^2\right)^{1/2}, \quad t\in(0,\,T).
\end{equation*}
The proof is completed, by using the embedding $H^1(0,\,T)\hookrightarrow C([0,\,T])$. \hfill $\Box$

\bigskip
We are now in a position to show the main result.
\begin{theorem}\label{thM_main}
Under the condition
\begin{equation}\label{M_cond_1}
\varepsilon \max_{1\leq i \leq M}\sum_{j\not=i} d_{ij}^{-2} < 1,
\end{equation}
which means that $1-2\beta-s/3\geq 0$, we have the following asymptotic expansion:
\begin{equation}\label{M_main}
u^s(x,\,t) = \sum_{j=1}^M \frac{C_j \alpha_j(t-c_0^{-1}|x-z_j|)}{4\pi |x-z_j|} + O\left(\varepsilon^{2-s}\right) + O\left(\varepsilon^{3-2s}\right)+ O\left(\varepsilon^{3-2\beta-s}\right) \quad \mathrm{as}\;\varepsilon\to 0
\end{equation}
for $x\in \mathbb R^3\setminus\overline D$ and $t\in (0,\,T)$, where the constants $C_j$'s are defined by \eqref{C_i} and $\{\alpha_j \}_{j=1}^M\subset H_0^1(0,\,T)$ is the unique solution of the linear algebraic system
\begin{equation}\label{Mq_i}
\alpha_i(t) + \sum_{\substack{j=1\\j\not= i}}^M \frac{C_j\alpha_j(t-c_0^{-1}|z_i-z_j|)}{4\pi |z_i-z_j|} = - \frac{\lambda(t-c_0^{-1}|z_i-z^*|)}{4\pi |z_i-z^*|}, \quad i=1,\,2,\,\cdots,\,M.
\end{equation}
\end{theorem}

{\bf Proof.} By \eqref{Mv_i_1} and \eqref{Mq_i}, we obtain from Lemma \ref{Mlinear} that
\begin{eqnarray*}\label{Mq_i1}
\left( \sum_{i=1}^M |\alpha_i(t) - C_i^{-1}v_i(t)|^2 \right)^{1/2} &\lesssim & \left( 1 - C \max_{1\leq i \leq M} \sum_{j\not=i} \frac{1}{4\pi |z_i-z_j|}\right)^{-1}\, \left( \sum_{i=1}^M \left(E_i^{(4)}\right)^2 \right)^{1/2} \nonumber\\
& = & O\left( M^{1/2}\,\max_{1\leq i \leq M}E_i^{(4)} \right) \nonumber\\
& = & O\left( M^{1/2}\,\left[\varepsilon + \varepsilon^2\,\left(O(d^{-2})+ O(M)\right)\right] \right),
\end{eqnarray*}
due to the fact in \cite{A-C-C-S-19} that $\max_{1\leq i \leq M} \sum_{j\not = i} d_{ij}^{-2} = O\left( d^{-2}\right) + O \left( M \right).$ For $x\in \mathbb R^3\setminus\overline D$ and $t\in (0,\,T)$, we have
\begin{eqnarray*}\label{M_main_1}
u^s(x,\,t) &=& \sum_{j=1}^M \int_{\partial D_j} \frac{\varphi_j(y,\,t-c_0^{-1}|x-y|)}{4\pi |x-y|}\,ds(y) \nonumber\\
&=& \sum_{j=1}^M \frac{1}{4\pi |x-z_j|}\int_{\partial D_j} \varphi_j(y,\,t-c_0^{-1}|x-z_j|)\,ds(y) \nonumber\\
&& \qquad +  \sum_{j=1}^M \int_{\partial D_j} \frac{1}{4\pi |x-y|}\left[\varphi_j(y,\,t-c_0^{-1}|x-y|) - \varphi_j(y,\,t-c_0^{-1}|x-z_j|)\right]\,ds(y) \nonumber\\
&& \qquad +  \sum_{j=1}^M \int_{\partial D_j} \varphi_j(y,\,t-c_0^{-1}|x-z_j|)\,\left[ \frac{1}{4\pi |x-y|} - \frac{1}{4\pi |x-z_j|}\right]\,ds(y) \nonumber\\
&=& \sum_{j=1}^M \frac{C_j \alpha_j(t-c_0^{-1}|x-z_j|)}{4\pi |x-z_j|} + O\left(C M\left[\varepsilon + \varepsilon^2\,\left(O(d^{-2})+O(M) \right)\right]\right) + O\left(M\varepsilon^2\right)\nonumber\\
&=& \sum_{j=1}^M \frac{C_j \alpha_j(t-c_0^{-1}|x-z_j|)}{4\pi |x-z_j|} + O\left(\varepsilon^{2-s}\right) + O\left(\varepsilon^{3-2s}\right)+ O\left(\varepsilon^{3-2\beta-s}\right).
\end{eqnarray*}
The proof is now complete. \hfill $\Box$

\section{Application to the effective medium theory}\label{effective-medium}
\setcounter{equation}{0}

In this section, we prove Theorem \ref{Main2} by utilizing the asymptotic expansion in the limit case that the holes are densely distributed and occupy a bounded domain.

\begin{figure}[htp]
	\centering
	\includegraphics[width=0.45\linewidth,height=4.5cm]{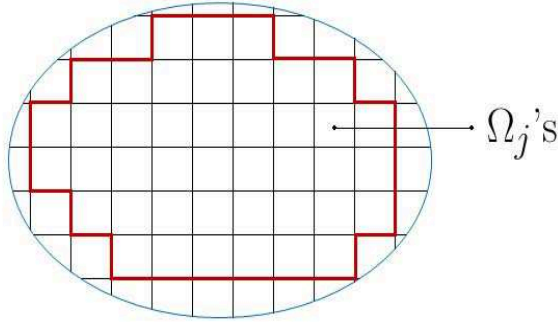}
	\caption{Bold red line encloses $\cup^{[a^{-1}]}_{j=1}\Omega_j$.}
	\label{omega_division}
\end{figure}

Let $\Omega$ be a bounded domain containing the holes $D_j,\,j=1,\,2,\,\cdots,\,M$. We divide $\Omega$ into $[a^{-1}]$ periodically subdomains $\Omega_j,\,j=1,\,2,\,\cdots,\,[a^{-1}]$ such that $\Omega_j$'s are disjoint and each $\Omega_j$ contains one single hole $D_j$ and has a volume $a$; see Figure \ref{omega_division}. We also assume that the holes $D_j,\,j=1,\,2,\,\cdots,\,M$ have the same shape. This means that $C_i=C_j$ for $i,\,j=1,\,2,\,\cdots,\,M$. Define $C:=C_j=\overline C\; a,$ where $\overline C$ is the scaled value of $C_j$.

Since $\Omega$ can have an arbitrary shape, the set of the cubes intersecting $\partial \Omega$ is not empty (unless if $\Omega$ has a simple shape as a cube). Later in our analysis, we will need the estimate of the volume of this set. Since each $\Omega_j$ has volume of the order $a$, its maximum radius is of the order $a^{\frac{1}{3}}$, and then the intersecting surfaces with $\partial \Omega$ has an area of the order $a^{\frac{2}{3}}$. As the area of $\partial \Omega$ is of the order one, we conclude that the number of such cubes will not exceed the order $a^{-\frac{2}{3}}$. Hence the volume of this set will not exceed the order $a^{-\frac{2}{3}}a=a^{\frac{1}{3}}$, as $a \rightarrow 0$. In particular $\mathrm{Vol}\,(\Omega\setminus\cup^{[a^{-1}]}_{j=1}\Omega_j)=O\left(a^{\frac{1}{3}}\right)$.

We consider the integral equation
\begin{equation}\label{med_1}
v(x,\,t) + \int_\Omega \overline C\, \frac{v(z,\,t-c_0^{-1}|x-z|)}{4\pi |x-z|} \,dz = -u^i(x,\,t), \quad (x,\,t)\in\Omega_T.
\end{equation}
Following the convolution quadrature based argument in \cite{Lubich1994}, we can prove that the equation \eqref{med_1} has a unique solution in $H_{0,\sigma}^r(0,\,T;\,L^2(\Omega))$ for $u^i\in H_{0,\sigma}^{r+2}(0,\,T;\,L^2(\Omega))$ with $r\in \mathbb Z_+$; see \cite{L-M2015}. As $u^i(x,\,t)$ given by \eqref{ui} with $z^*\not\in\overline\Omega$ is sufficiently smooth for $(x,\,t)\in \Omega_T$, we have $v\in H_{0,\sigma}^2(0,\,T;\,L^2(\Omega))$, and hence $ v\in C^1\left([0,\,T];\,L^2(\Omega) \right)$ by Sobolev embedding $H^2(0,\,T)\hookrightarrow C^1[0,\,T]$.

Define
\begin{equation}\label{def_V}
V(x,\,t):=
\begin{cases}
v(x,\,t) & \mathrm{in}\;\Omega_T,\\
-u^i(x,\,t) - \displaystyle\int_\Omega \overline C\, \frac{v(z,\,t-c_0^{-1}|x-z|)}{4\pi |x-z|} \,dz & \mathrm{in}\; (\mathbb R^3\setminus\overline \Omega)_T,
\end{cases}
\end{equation}
and set $W(x,\,t) := -u^i(x,\,t) - V(x,\,t)$. Then $W(x,\,t)$ satisfies
\begin{equation}\label{W-1}
\begin{cases}
(c_0^{-2}\partial_{tt} - \Delta + \overline C \chi_\Omega)W = - \overline C \chi_\Omega u^i(x,\,t) & \mathrm{in}\;\mathbb R^3\times (0,\,T), \\
W(x,\,0)=0,\,W_t(x,\,0)=0 & \mathrm{in}\; \mathbb R^3.
\end{cases}
\end{equation}
If we define $U:=W+u^i$, we also have
\begin{equation}\label{U}
\begin{cases}
(c_0^{-2}\partial_{tt} - \Delta + \overline C \chi_\Omega)U =0 & \mathrm{in}\;\mathbb R^3\times (0,\,T), \\
U(x,\,0)=0,\,U_t(x,\,0)=0 & \mathrm{in}\; \mathbb R^3.
\end{cases}
\end{equation}
The main result of this section is stated as follows.

\begin{theorem}\label{th_eff}
For any fixed $x\in\mathbb R^3\setminus\overline \Omega$ and $t\in (0,\,T)$, we have the estimate
\begin{equation}\label{eff}
W(x,\,t) = u^s(x,\,t) + O\left( a^{\frac{1}{3}} \right)  \quad \mathrm{as}\,a\to 0,
\end{equation}
or equivalently,
\begin{equation}\label{eff1}
U(x,\,t) = u(x,\,t) + O\left( a^{\frac{1}{3}} \right)  \quad \mathrm{as}\,a\to 0
\end{equation}
where $u^s(x,\,t)$ is the solution to \eqref{ibvp} and $u=u^i+u^s$.
\end{theorem}

{\bf Proof.} First, we show the regularity of $v(x,\,t)$. From \eqref{med_1}, we see
$$
\vert v(x,\,t)\vert \lesssim \int_\Omega |x-z|^{-1} \,|v(z, t-c_0^{-1}|x-z|) | \,dz + |u^i(x,\,t)|,
$$
and hence
$$
\vert v(x,\,t)\vert \lesssim \left(\int_\Omega |x-z|^{-2}\,dz\right)^\frac{1}{2} \Vert v \Vert_{C\left([0,\,T];\,L^2(\Omega) \right)}+O(1).
$$
So $v$ is in $C\left([0,\,T]; \,L^\infty(\Omega) \right)$. In addition, by taking the derivative on the both sides of \eqref{med_1}, we can also get
\begin{eqnarray*}\label{med_4}
|\partial_{x_j}v(x,\,t) | &\lesssim & \int_\Omega \frac{|v(z,\,t-c_0^{-1}|x-z|)|}{4\pi |x-z|^2} \,dz + \int_\Omega \frac{ |\partial_t v(z,\,t-c_0^{-1}|x-z|)|}{4\pi |x-z|} \,dz + |\partial_t u^i(x,\,t)|  \nonumber\\
&\lesssim & \int_\Omega |x-z|^{-2}\,dz\, \Vert v \Vert_{C\left([0,\,T];\,L^\infty(\Omega)\right)} + \left(\int_\Omega |x-z|^{-2}\,dz\right)^\frac{1}{2} \Vert \partial_t v \Vert_{C\left([0,\,T];\,L^2(\Omega) \right)}+O(1) \nonumber\\
&=& O(1).
\end{eqnarray*}
This means that $\partial_{x_j} v\in C\left([0,\,T];\,L^\infty(\Omega) \right)$. So we obtain that $v\in C\left([0,\,T];\, W^{1,\infty}(\Omega)\right)$. Analogously, we can also prove that $\partial_t v\in C\left([0,\,T];\, W^{1,\infty}(\Omega)\right)$

\medskip
Next, we estimate $\sum_{j=1}^M|\alpha_j(t)-v(z_j,\,t)|^2$, where $\{ \alpha_j(t) \}_{j=1}^M$ is the solution to the linear algebraic system \eqref{Mq_i}. To this end, we rewrite the integral equation \eqref{med_1} at $x=z_l$ for $1\leq l \leq M$ as
\begin{equation}\label{med_1_2}
v(z_l,\,t) + \sum_{\substack{j=1\\j\neq l}}^{M} \overline C\, {a} \frac{v(z_j,\,t-c_0^{-1}|z_l-z_j|)}{4\pi |z_l-z_j|}  = - u^i(z_l,\,t) + \mathcal A + \mathcal A_l + \mathcal B_l,
\end{equation}
where
\begin{eqnarray*}
\mathcal A &:= & - \int_{\Omega \setminus\left(\cup^{[a^{-1}]}_{j=1}\Omega_j\right)} \overline C\,\frac{v(z,\,t-c_0^{-1}|z_l-z|)}{4\pi |z_l-z|} \,dz,\\
\mathcal A_l & := & - \int_{\Omega_l} \overline C\,\frac{v(z,\,t-c_0^{-1}|z_l-z|)}{4\pi |z_l-z|}\, dz,\\
\mathcal B_l & := & - \sum_{\substack{j=1\\j\neq l}}^{[a^{-1}]}\overline C  \int_{\Omega_j} \,\frac{v(z,\,t-c_0^{-1}|z_l-z|)}{4\pi |z_l-z|} \,dz +  \sum_{\substack{j=1\\j\neq l}}^{[a^{-1}]}\, \overline C\, a\, \frac{v(z_j,\,t-c_0^{-1}|z_l-z_j|)}{4\pi |z_l-z_j|}.
\end{eqnarray*}

Since $v\in C\left([0,\,T];\,L^\infty(\Omega) \right)$, we have
$$\mathcal A_l  = O\left( \int_{\Omega_l} \vert z-z_l\vert^{-1}\, dz\right),$$ and hence, by a scaling, we get the estimate
\begin{equation}\label{est_A_l}
\mathcal A_l =O\left(a^{\frac{2}{3}}\right)\quad \mathrm{as}\; \varepsilon \ll 1.
\end{equation}
Let us estimate $\mathcal B_l$. As $\vert \Omega_l\vert= a $, we have
$$
\mathcal B_l=- \sum_{\substack{j=1\\j\neq l}}^{[a^{-1}]} \overline C \int_{\Omega_j} \left[\frac{v(z,\,t-c_0^{-1}|z_l-z|)}{4\pi |z_l-z|}  - \frac{v(z_j,\,t-c_0^{-1}|z_l-z_j|)}{4\pi |z_l-z_j|} \right] \,dz.
$$
Write the above integrand as
\begin{eqnarray*}
&&\frac{v(z,\,t-c_0^{-1}|z_l-z|)}{4\pi |z_l-z|}  - \frac{v(z_j,\,t-c_0^{-1}|z_l-z_j|)}{4\pi |z_l-z_j|} \\
&=& v(z,\,t-c_0^{-1}|z_l-z|)\left[\frac{1}{4\pi |z_l-z|}  - \frac{1}{4\pi |z_l-z_j|} \right] \\
& & + \frac{1}{4\pi |z_l-z_j|} \left[v(z,\,t-c_0^{-1}|z_l-z|) - v(z_j,\,t-c_0^{-1}|z_l-z_j|) \right].
\end{eqnarray*}
Then we see
\begin{equation}\label{est_B_l}
\mathcal B_l=O\Big(\sum_{\substack{j=1\\j\neq l}}^{[a^{-1}]} d_{lj}^{-2}\Big)\, a^{\frac{4}{3}}=O\left(a^{\frac{1}{3}}\right).
\end{equation}

Let us estimate the term $\mathcal A$. We distinguish the following two cases:
\begin{itemize}

\item[(1)] The point $z_l$ is away from the boundary $\partial \Omega$ and so $|z_l-z|^{-1} $ is bounded in $z$ near the boundary. In this case,
we have $\mathcal A =O\left(\mathrm{Vol}\,(\Omega\setminus\cup^{[a^{-1}]}_{j=1}\Omega_j)\right)=O\left(a^{\frac{1}{3}}\right)$.

\item[(2)] The point $z_l$ is located near one of the $\Omega_j$'s touching the boundary $\partial \Omega $. In this case, we split the estimate into two parts. By $ N_l$ we denote the part that involves 
$\Omega_j $'s close to $z_l $, and we denote the remaining part by $F_l $. The integral over $F_l $ can be estimated in a manner similar to the case $(1)$ discussed above. 
Also note that $F_l \subset \Omega\setminus \cup_{j=1}^{[a^{-1}]} \Omega_j $ and so $\mathrm{Vol}\,(F_l) $ is of the order $a^{\frac{1}{3}} $ as $a \rightarrow 0 $.

\medskip
To estimate the integral over $N_l $, we first estimate the number of $\Omega_j$'s close to $z_l$. We observe that the $\Omega_j $'s close to $z_l$ are located near a small region of the boundary $\partial \Omega $. 
Since we assume that the boundary is smooth enough, this region can be assumed to be flat and centered at $z_l$.  We now divide this flat region into concentric squared layers (centered at $z_l$); see Figure \ref{layers}. 
Observe that as this flat region is of order $1$, in term of the parameter $a$, and the maximum radius of the squares (or the $\Omega_j$'s) is $a^{\frac{1}{3}}$, 
then the number of the layers is at most of the order $[a^{-\frac{1}{3}}]$.
In this case, we have at most $(2n+1)^{2} $ squares (and hence cubes intersecting the surface) in the $n$ first layers, for $n=0,\,\dots,\,[a^{-\frac{1}{3}}] $. 
So the number of holes in the $n^{th} $ layer $(n \neq 0)$ will be at most $[(2n+1)^{2}-(2n-1)^{2}] $ and their distance from $\Omega_l$ is at least $n\left(a^{\frac{1}{3}}-\frac{a}{2} \right)$.
\begin{figure}[htp]
	\centering
	\includegraphics[width=0.4\linewidth,height=4.5cm]{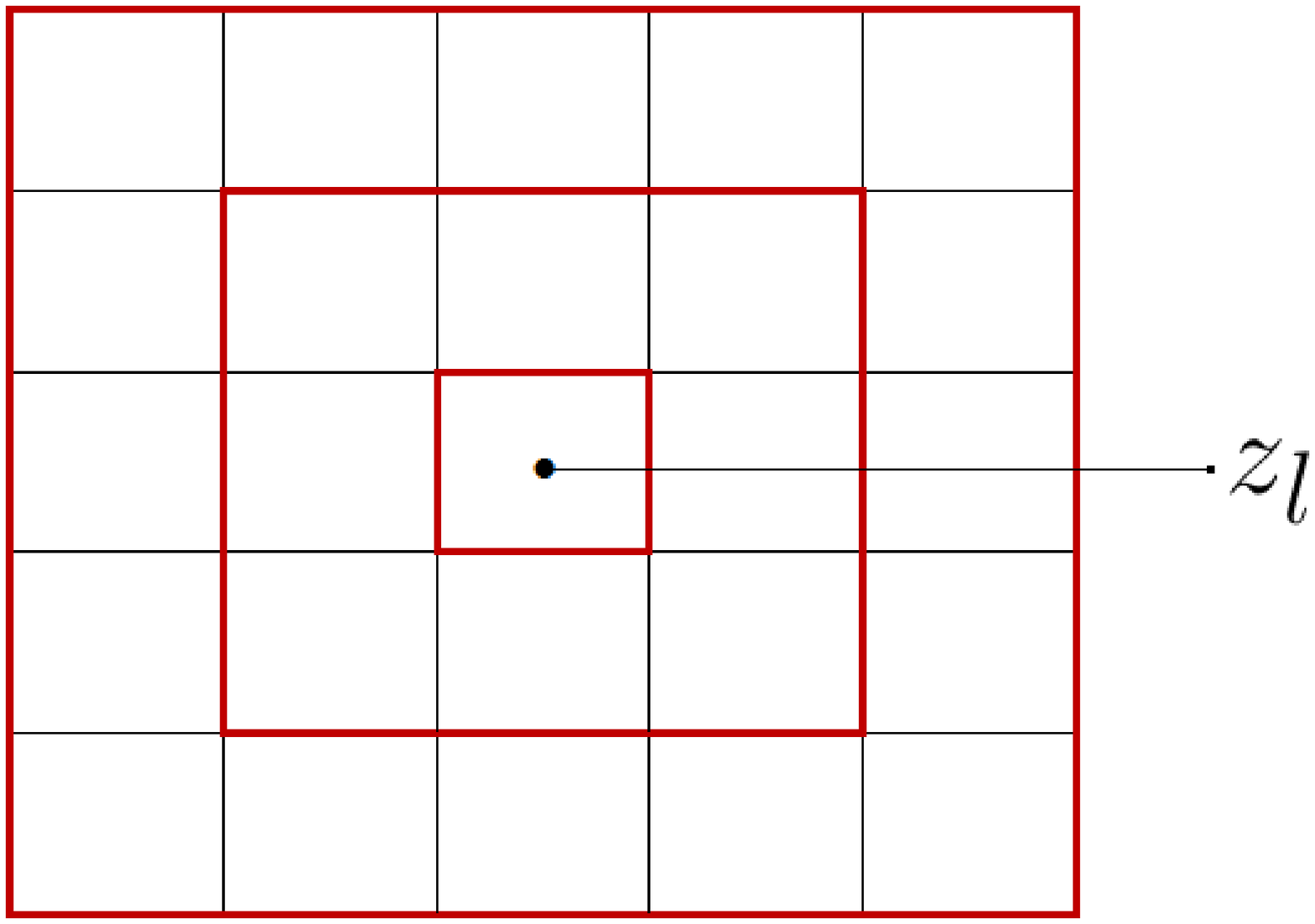}
	\caption{Concentric squared layers centered at $z_l$ (with three layers, i.e. $n=2$).}
	\label{layers}
\end{figure}
\end{itemize}

Therefore we can write
\begin{eqnarray*}
\vert \mathcal A \vert &=& \left\vert \int_{\Omega \setminus \left(\cup_{j=1 }^{[a^{-1}]} \Omega_j\right)} \overline C\, \frac{v(z,\,t-c_0^{-1}|z_l-z|)}{4\pi |z_l-z|} \, dz \right\vert\\
&\leq& \left\vert \int_{N_l} \overline C\, \frac{v(z,\,t-c_0^{-1}|z_l-z|)}{4\pi |z_l-z|}\, dz \right\vert + \left\vert\int_{F_l} \overline C\, \frac{v(z,\,t-c_0^{-1}|z_l-z|)}{4\pi |z_l-z|} \, dz \right\vert \\
&\leq& \sum_{m=1}^{[a^{-\frac{1}{3}}]}\overline C\,  \Vert v\Vert_{C\left([0,\,T];\,L^\infty(\Omega) \right)} \,\mathrm{Vol}\,(\Omega_m)\, \frac{1}{d_{ml}} +  
\overline C\, \Vert v \Vert_{C\left([0,\,T];\,L^\infty(\Omega) \right)}\, \mathrm{Vol}\,(F_l) \\
&\leq& O\Big(a \sum_{m=1}^{[a^{-\frac{1}{3}}]} \frac{1}{d_{ml}} + \overline C\, a^{\frac{1}{3}}\Big) \\
&\leq& O\Big( a \Big[(2n+1)^2-(2n-1)^2 \Big]\, \frac{1}{n \Big( a^{\frac{1}{3}}-\frac{a}{2} \Big)} + \overline C\,a^{\frac{1}{3}}\Big)\\
&=& O\Big( a \,O(a^{-\frac{2}{3}}) + O(a^{\frac{1}{3}})\Big),
\end{eqnarray*}
and hence
\begin{equation}\label{est-vol-2a}
\vert \mathcal A \vert = O\left(a^{\frac{1}{3}}\right).
\end{equation}

Gathering the estimates \eqref{est_A_l}, \eqref{est_B_l} and \eqref{est-vol-2a}, we have
$$
\sum_l \left(\vert \mathcal A \vert^2 +\vert \mathcal A_l \vert^2 +\vert \mathcal B_l\vert^2 \right)=
O\left(Ma^{\frac{2}{3}}+ Ma^{\frac{4}{3}}\right)=O\left(a^{-\frac{1}{3}}\right).
$$
Using the invertibility property and the estimate \eqref{Mlinear_2} for the algebraic system \eqref{Mlinear_1}, we deduce the following estimate:
\begin{equation}\label{med_8}
\sum_{j=1}^M |\alpha_j(t) - v(z_j,\,t) |^2 = O(a^{-\frac{1}{3}}) \quad \mathrm{as}\; a \to 0.
\end{equation}

Finally, we estimate $|W(x,\,t) - u(x,\,t)|$. Let $x$ be away from $\Omega \cup \{z^*\}$. Recall that
\begin{equation*}
W(x,\,t) = -u^i(x,\,t) - V(x,\,t) = \int_\Omega \overline C\, \frac{v(z,\,t-c_0^{-1}|x-z|)}{4\pi |x-z|}\,dz,
\end{equation*}
and rewrite it as
\begin{equation*}
W(x,\,t) = \sum_{j=1}^{[a^{-1}]} |\Omega_j| \, \overline C \,\frac{v(z_j,\,t-c_0^{-1}|x-z_j|)}{4\pi |x-z_j|} + \mathcal D
\end{equation*}
with
\begin{equation*}
\mathcal D = \sum_{j=1}^{[a^{-1}]} \int_{\Omega_j}\overline C\,\left[\frac{v(z,\,t-c_0^{-1}|x-z|)}{4\pi |x-z|} -  \frac{v(z_j,\,t-c_0^{-1}|x-z_j|)}{4\pi |x-z_j|} \right]\,dz + \int_{\Omega \setminus\left(\cup^{[a^{-1}]}_{j=1}\Omega_j\right)} \overline C\,\frac{v(z,\,t-c_0^{-1}|x-z|)}{4\pi |x-z|} \,dz.
\end{equation*}
Following the similar steps as for estimating $\mathcal B_l$ and $\mathcal A$, and as the integrands are smooth here, it can be easily proved that $\mathcal D = O(a^{\frac{1}{3}})$ as $a \to 0$. Then we have
\begin{equation*}
W(x,\,t) = \sum_{j=1}^{[a^{-1}]} \overline C \, |\Omega_j| \frac{\alpha_j(t-c_0^{-1}|x-z_j|)}{4\pi |x-z_j|} + \mathcal E + O(a^{\frac{1}{3}})
\end{equation*}
with
\begin{equation*}
\mathcal E := -\sum_{j=1}^{[a^{-1}]} \overline C \, |\Omega_j| \,\frac{\alpha_j(t-c_0^{-1}|x-z_j|) - v(z_j,\,t-c_0^{-1}|x-z_j|)}{4\pi |x-z_j|}.
\end{equation*}
The term $\mathcal E$ can be estimated as
\begin{equation*}
\mathcal E = O\left(a\, M^{1/2}a^{-1/6}\right)=O\left(a^{\frac{1}{3}}\right).
\end{equation*}
Hence, we conclude from Theorem \ref{thM_main} that
\begin{equation}\label{med_final}
W(x,\,t) = u^s(x,\,t) + O\left(a^{\frac{1}{3}}\right) \quad \mathrm{as}\; a \to 0.
\end{equation}
The proof is now complete. \hfill $\Box$

\section{Numerical examples} \label{numer}
\setcounter{equation}{0}

In this section, we show three numerical examples to verify our theoretical results in Theorems \ref{Main} and \ref{Main2}. Examples \ref{exam_1} and \ref{exam_2} are presented to illustrate the effectiveness of the asymptotic expansion \eqref{M1M_main}, while Example \ref{exam_4} is devoted to testing the approximation \eqref{Meff}. To numerically solve the scattering problem \eqref{ibvp}, we truncate the infinite domain $\mathbb R^3\setminus\overline D$ by a large enough spherical domain $\Omega_B$ such that $\overline D\subset \Omega_B$ and the scattered field on $\partial \Omega_B$ in a finite time interval $(0,\,T)$ is zero by Huygens' principle. That is, we consider the following initial boundary value problem in a bounded domain:
\begin{equation}\label{ibvp_bd}
\begin{cases}
c_0^{-2}u^s_{tt} - \Delta u^s=0 & \textrm{ in } (\Omega_B\setminus\overline D)_T, \\
u^s=0 & \textrm{ on } (\partial \Omega_B)_T, \\
u^s=-u^i & \textrm{ on } (\partial D)_T, \\
u^s|_{t=0}=0,\;u^s_t|_{t=0}=0 & \textrm{ in } \Omega_B\setminus\overline D.
\end{cases}
\end{equation}

In all numerical examples, we take the causal signal $\lambda(t)$ in the incident wave \eqref{ui} as 
\begin{equation*}
\lambda(t) = 
\begin{cases}
\exp\left(-t^{-2}\right), & t>0, \\
0, & t\leq 0
\end{cases}
\end{equation*}
and set the wave speed $c_0$ of the background medium as $c_0=1$. The domain $\Omega_B$ is fixed as a ball of radius $R = 1.2$ centered at the origin. The holes $D_j,\,j=1,\,\cdots,\,M$ are balls of radius $\varepsilon$ with different centers, and then the capacitance for each hole $D_j$ is $C_j=4\pi \varepsilon$.

\begin{example}\label{exam_1}
Let $D$ be a small spherical hole with the radius $\varepsilon$ and center at $(0.1,\,0,\,0)$. Set $z^* = (0.15,\,0,\,0)$ and $T = 1$.
\end{example}

We solve the scattering problem \eqref{ibvp_bd} by using the finite element method and take its numerical solution as the exact one. Let $\Gamma_r$ be the sphere of radius $r$ centered at the origin, that is, $\Gamma_r=\{x\in\mathbb R^3:\,|x|=r\}$, where we compare the numerical solution of \eqref{ibvp_bd} via the finite element method (FEM) with the asymptotic approximation computed by \eqref{M1M_main}. To numerically verify the effectiveness of the asymptotic expansion \eqref{M1M_main} and the convergence of the asymptotic approximation as $\varepsilon\to 0$, we test the cases of different radii $\varepsilon$. 

\begin{figure}[htp]
	\begin{center}
		\includegraphics[width=0.45\textwidth,height=4.5cm]{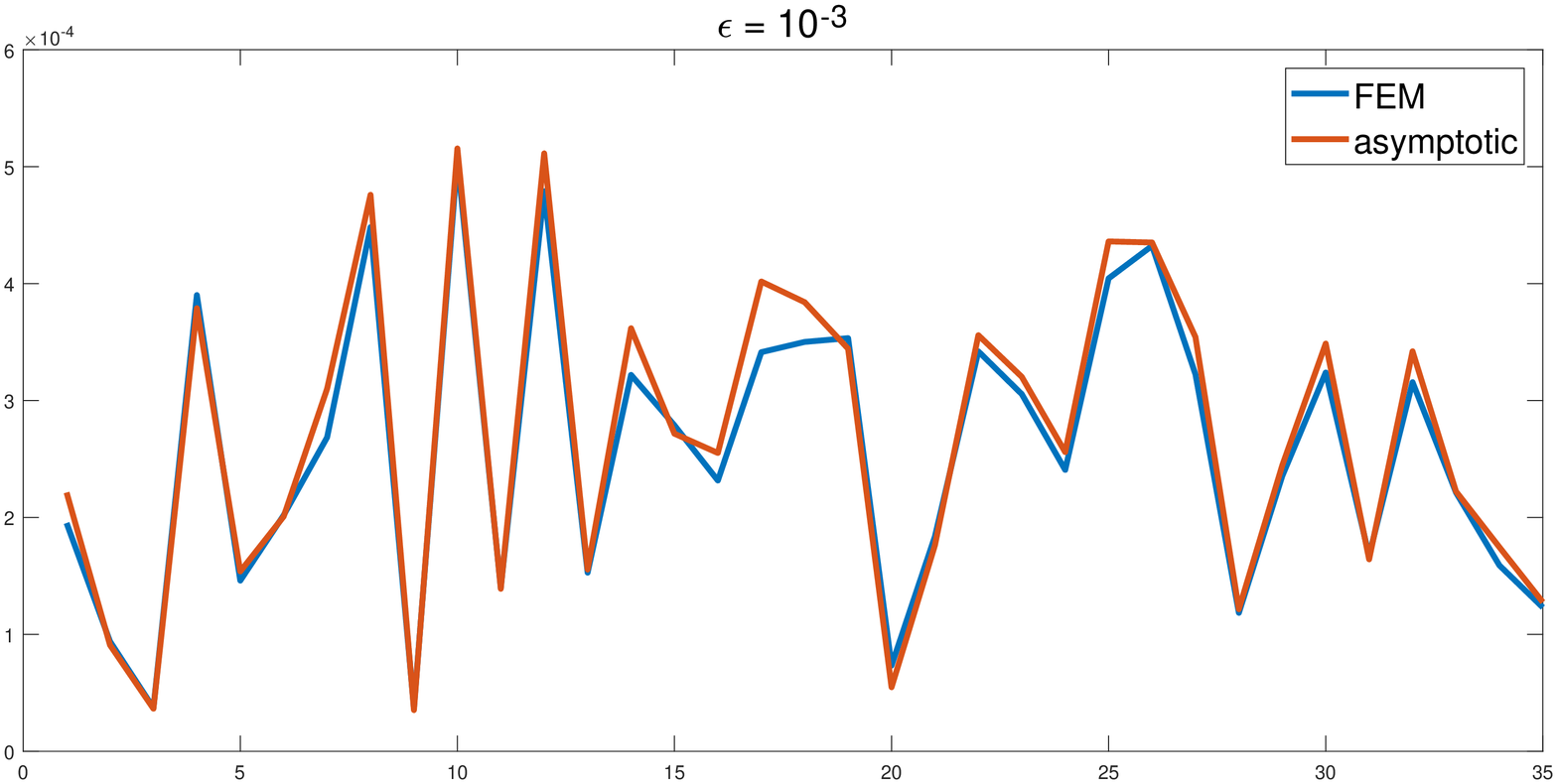}
		\includegraphics[width=0.45\textwidth,height=4.5cm]{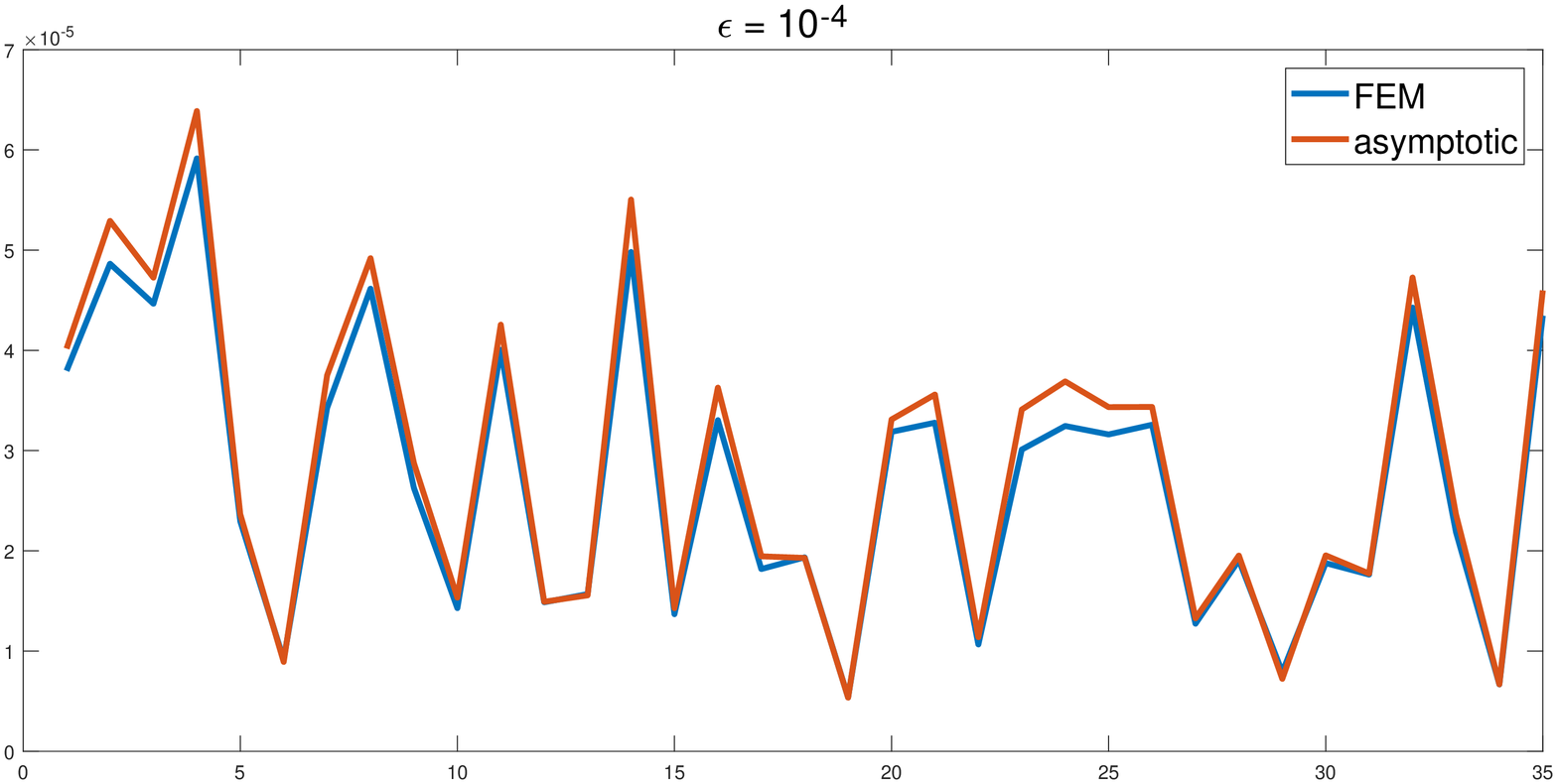}
		\includegraphics[width=0.45\textwidth,height=4.5cm]{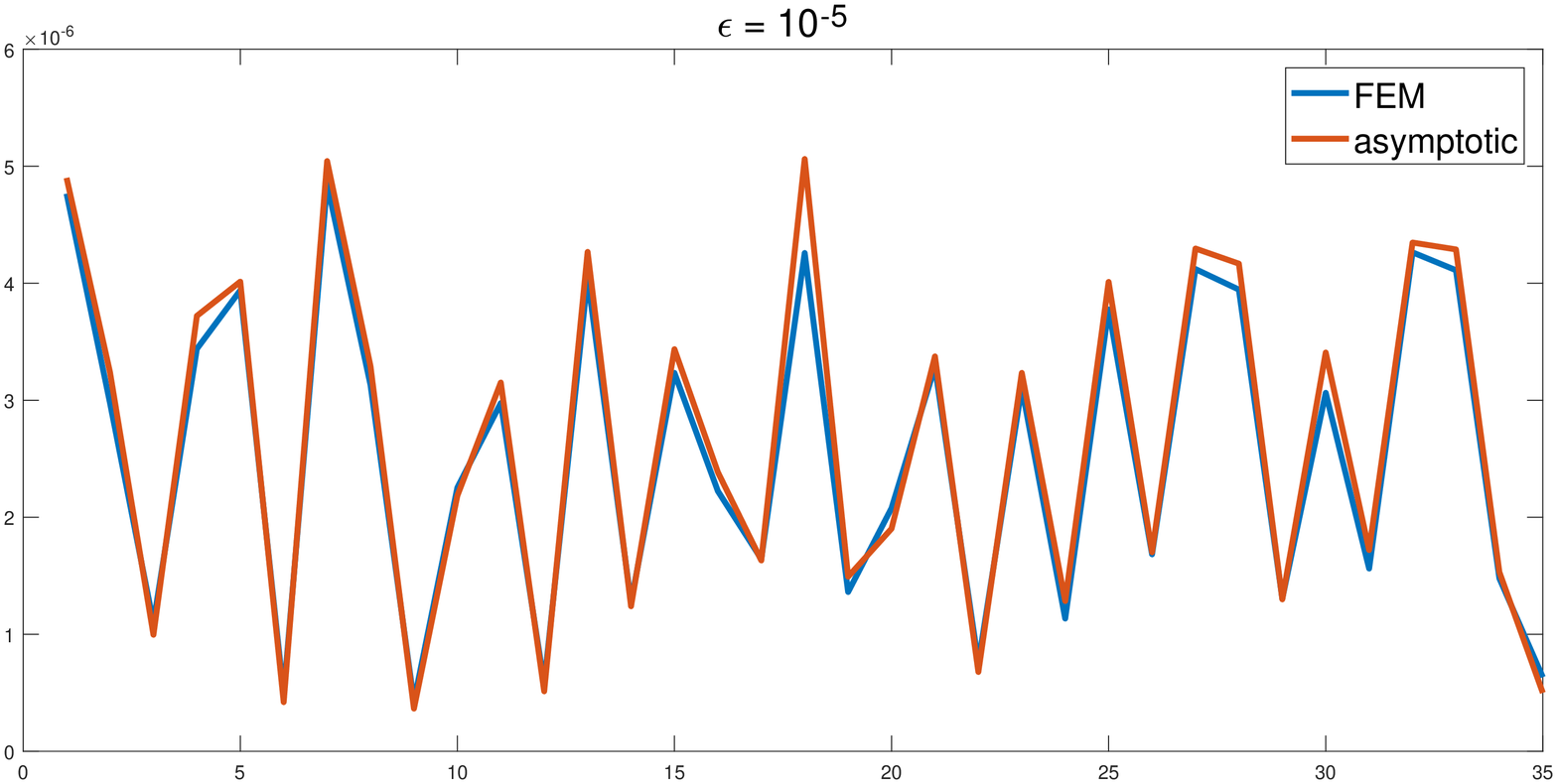}
		\includegraphics[width=0.45\textwidth,height=4.5cm]{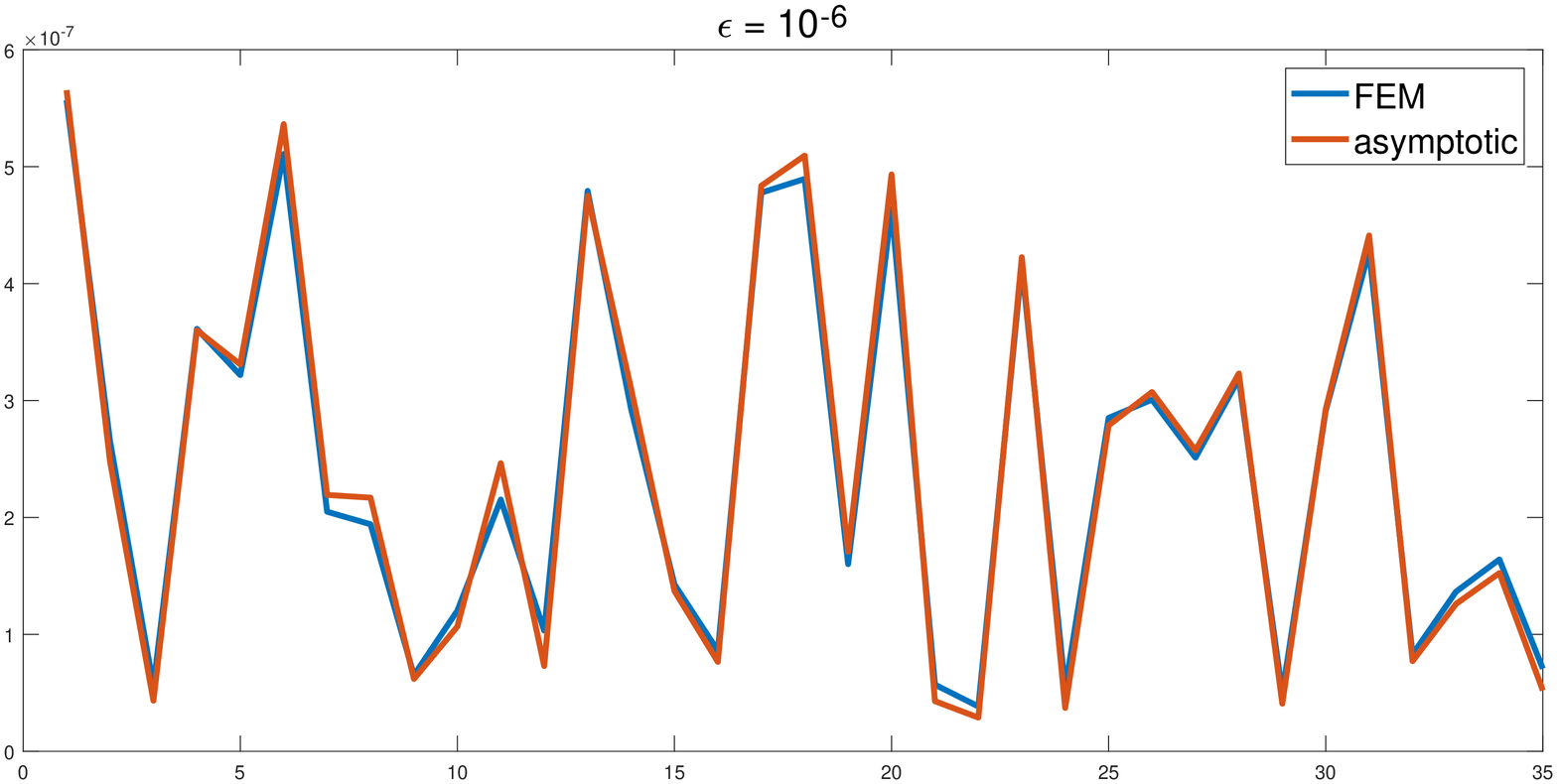}
		\caption{One hole case: the scattered field is computed on $\Gamma_r$ with $r=0.2$. The blue line denotes the numerical solution by FEM, and the red line stands for the asymptotic approximation.}
		\label{one}
	\end{center}
\end{figure}

\begin{figure}[htp]
	\begin{center}
		\includegraphics[width=0.45\textwidth,height=4.5cm]{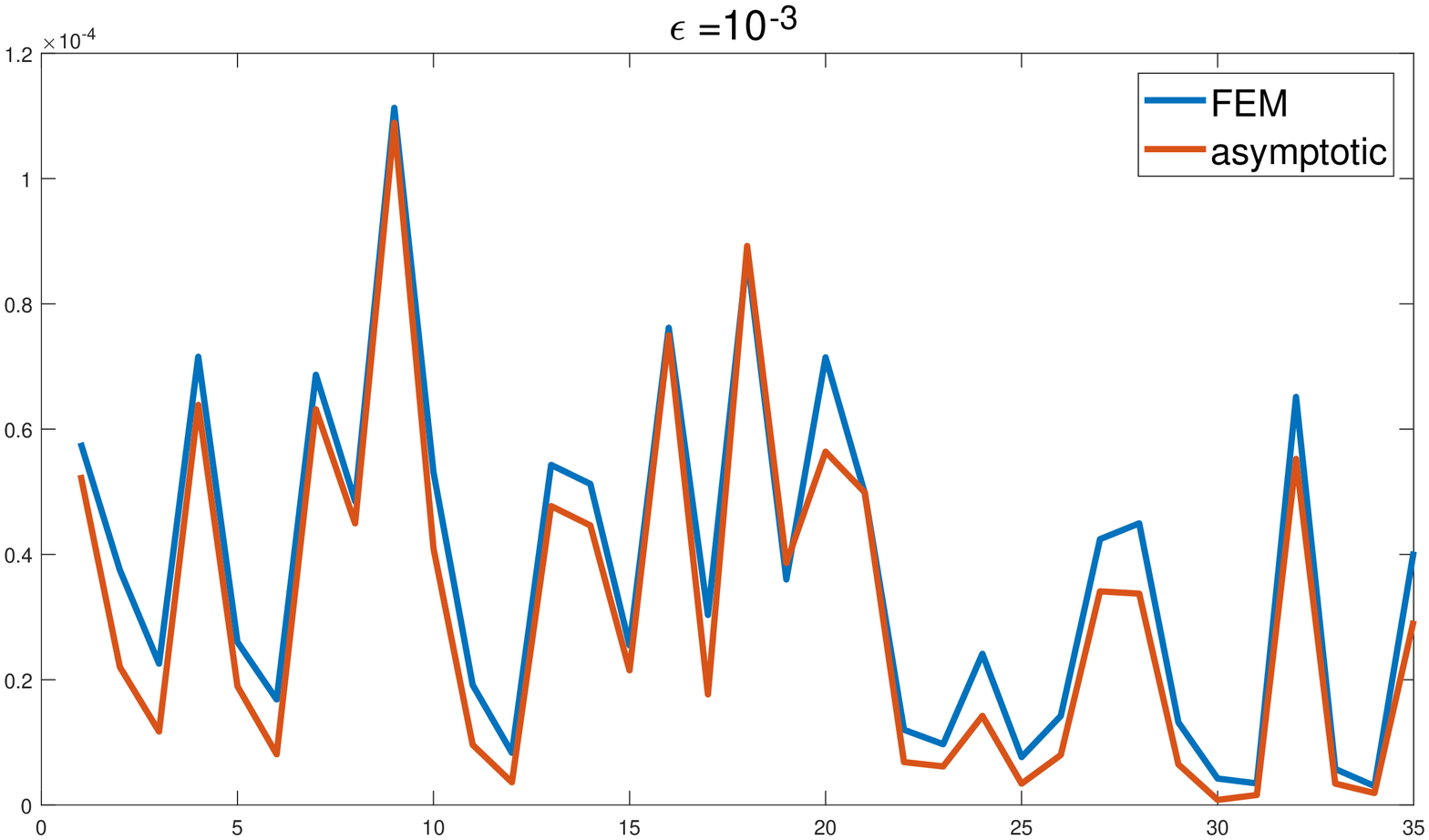}
		\includegraphics[width=0.45\textwidth,height=4.5cm]{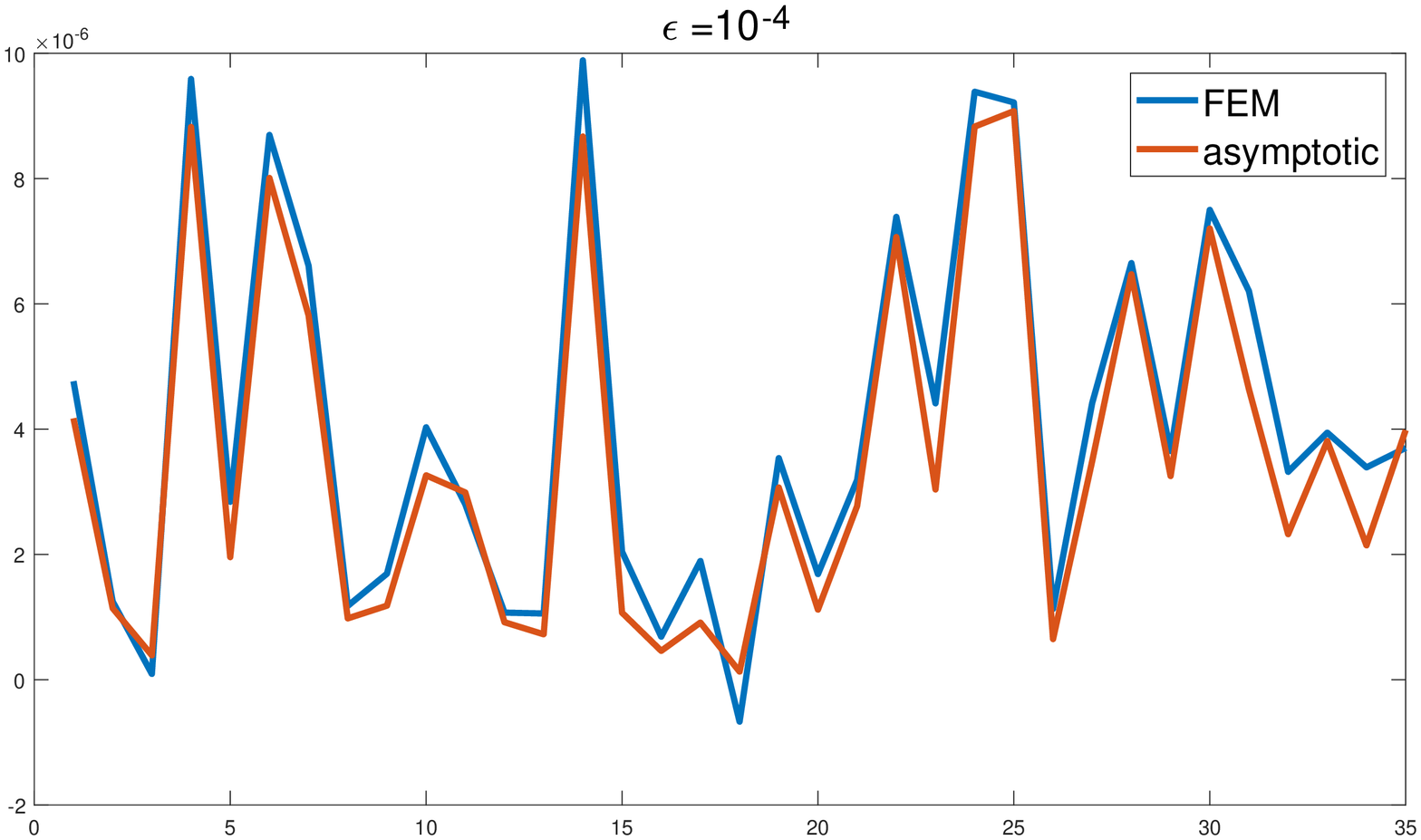}
		\includegraphics[width=0.45\textwidth,height=4.5cm]{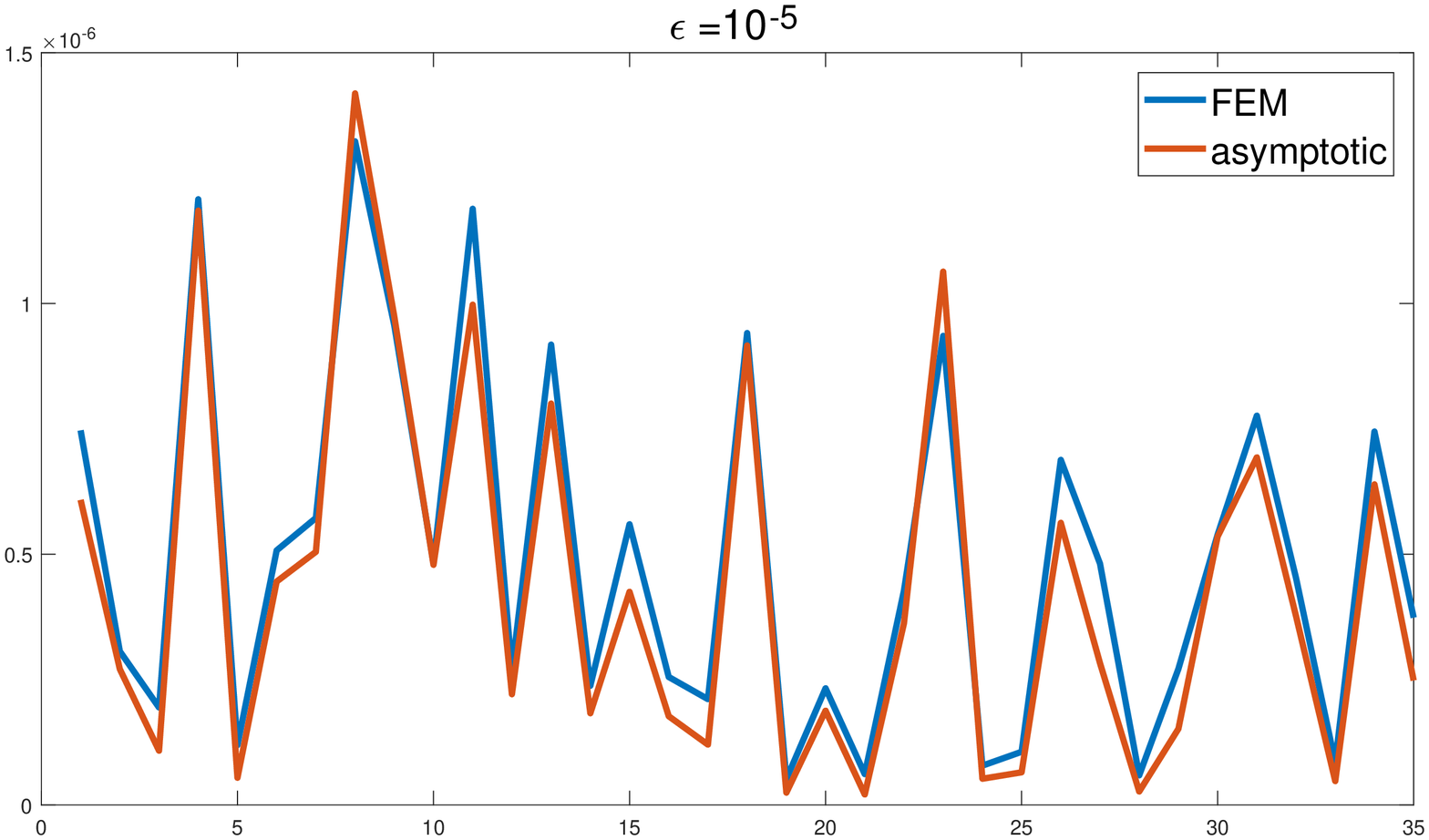}
		\includegraphics[width=0.45\textwidth,height=4.5cm]{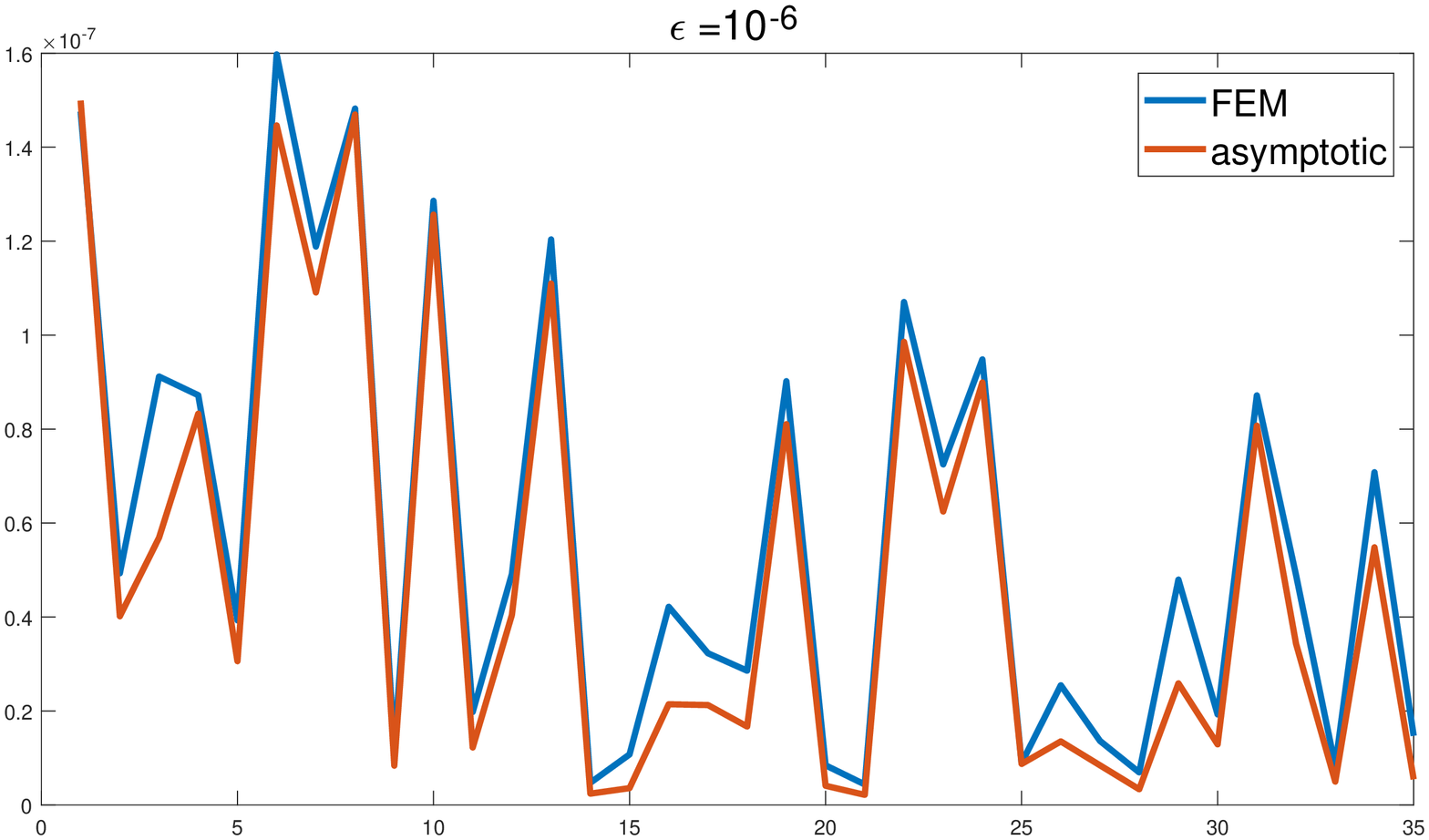}
		\caption{One hole case: the scattered field is computed on $\Gamma_r$ with $r=0.3$. The blue line denotes the numerical solution by FEM, and the red line stands for the asymptotic approximation.}
		\label{one2}
	\end{center}
\end{figure}

In Figure \ref{one}, we show the numerical results of the scattered field $u^s(x,\,t)$ on $\Gamma_r$ with $r=0.2$ for $\varepsilon=10^{-3},\,10^{-4},\,10^{-5},\,10^{-6}$. It can be easily observed that the dominant term of the asymptotic expansion gives a good approximation of the scattered field with reasonable errors and the approximation is evidently improved as the radius $\varepsilon$ becomes smaller. In Figure \ref{one2}, we also show the numerical results of the scattered field $u^s(x,\,t)$ on $\Gamma_r$ with $r=0.3$. We observe that the error of the asymptotic approximation becomes large as the observation points are away from the hole, which is reasonable from the derivation of the asymptotic expansion.

\begin{figure}[htp]
	\centering
	\includegraphics[width=0.4\linewidth,height=4.5cm]{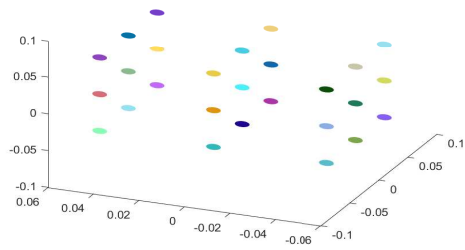}
	\caption{The distribution of $27$ small holes.}
	\label{fig:manyobstacles3d}
\end{figure}

\begin{example}\label{exam_2}
Let $D$ be the union of $27$ small spherical holes of radius $\varepsilon$ distributed in the cube $[-0.05,\,0.05]^3$; see Figure \ref{fig:manyobstacles3d} for the distribution of holes. Set $z^* = (0.123,\,0,\,0)$ and $T = 0.8$.
\end{example}

\begin{figure}[htp]
	\begin{center}
		\includegraphics[width=0.45\textwidth,height=4.5cm]{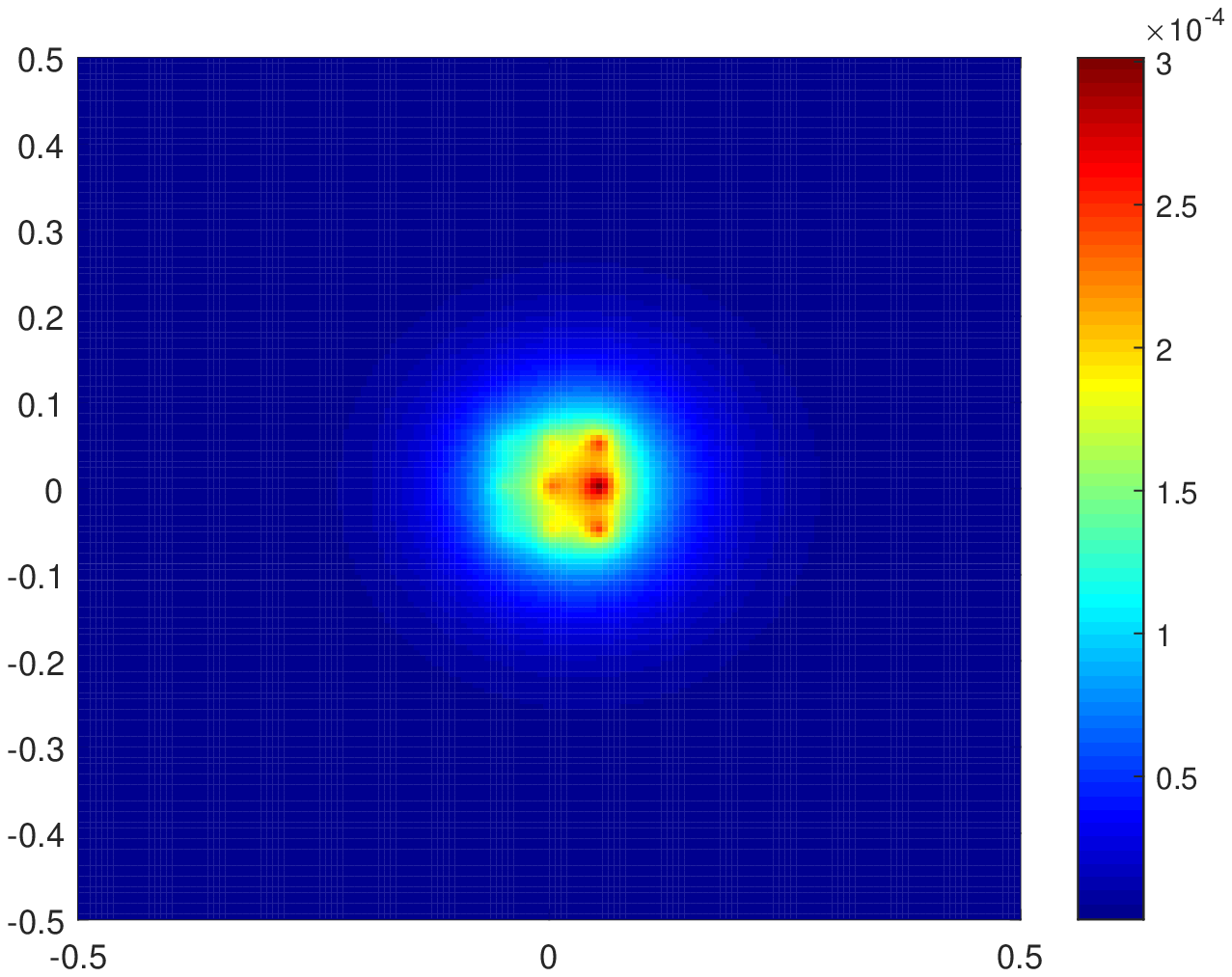}
		\includegraphics[width=0.45\textwidth,height=4.5cm]{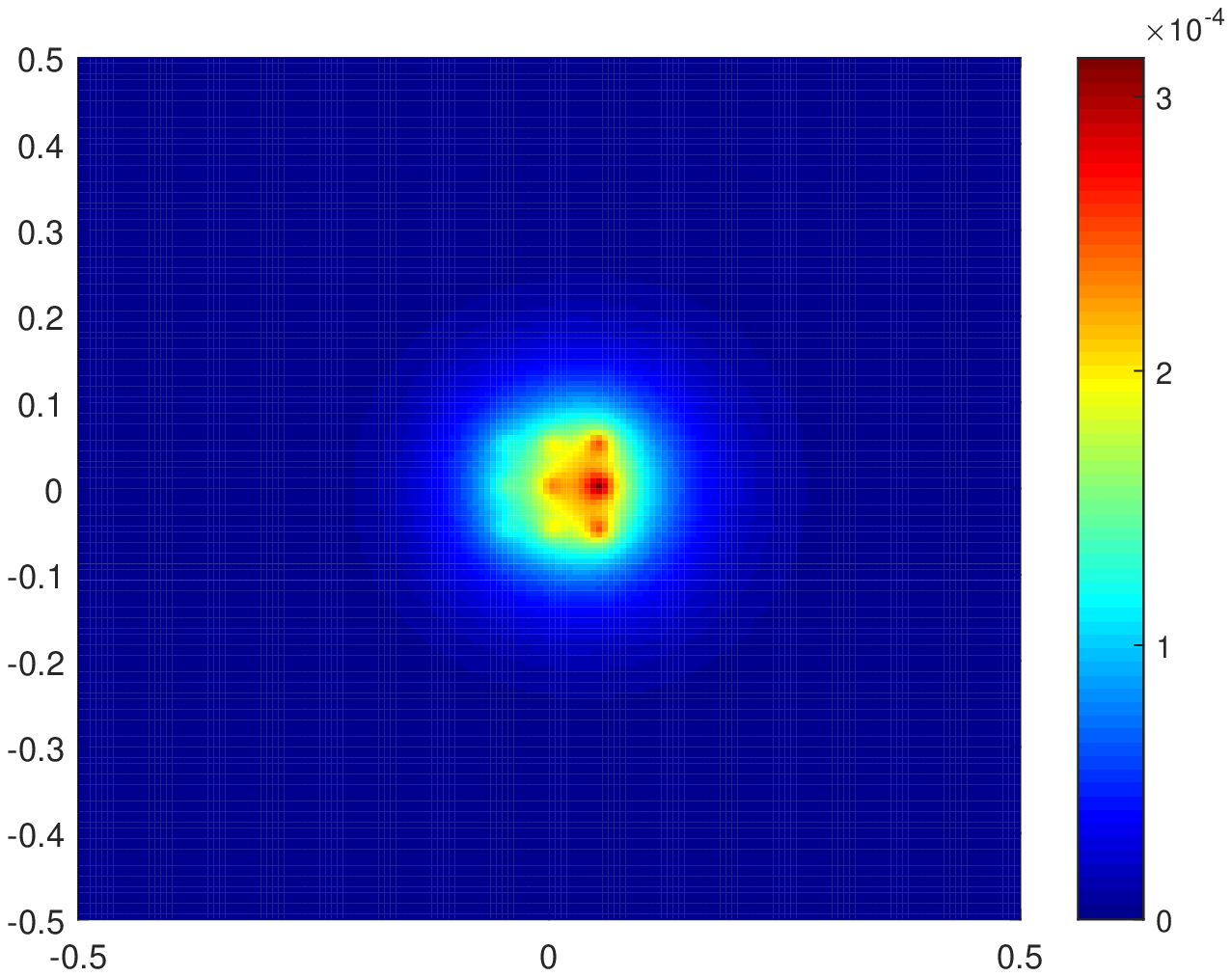}
		\caption{The case of $27$ holes: the scattered field is computed in the plane $z=0.01$. The left one shows the numerical solution by FEM, and the right one shows the numerical solution by the asymptotic approximation.}
		\label{many}
	\end{center}
\end{figure}

We test the case of radius $\varepsilon=10^{-5}$. The scattered field is computed in the planar domain $\{(x,\,y,\,z): -0.5\leq x,\,y \leq 0.5,\,z=0.01\}$, and we display in Figure \ref{many} the numerical solution obtained by solving \eqref{ibvp_bd} via the finite element method and the asymptotic approximation via \eqref{M1M_main}. We conclude from the above numerical results that the asymptotic expansion we derived could be used to approximately compute the scattered wave by a cluster of small holes. The complexity and computation time in using asymptotic approximation are much less than those for the finite element method.

\medskip
Finally, we show the performance of the approximation \eqref{Meff} by comparing the solutions to the original scattering problem \eqref{ibvp} and the effective medium problem \eqref{MW-1}.

\begin{figure}[htp]
	\centering
	\includegraphics[width=0.4\linewidth,height=4.5cm]{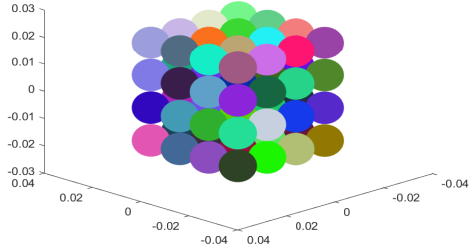}
	\caption{The distribution of $64$ small holes.}
	\label{fig:64obstacles_2}
\end{figure}

\begin{example}\label{exam_4}
Let $D$ be the union of $64$ small spherical holes of radius $\varepsilon=0.0055$, which are densely distributed in $\Omega=[-1.8\times 10^{-2},\,1.8\times 10^{-2}]^3$, see Figure \ref{fig:64obstacles_2}. Set $z^* = (0.0243,\,0,\,0)$ and $T = 0.8$.
\end{example}

\begin{figure}[htp]
	\begin{center}
		\includegraphics[width=0.45\textwidth,height=4.5cm]{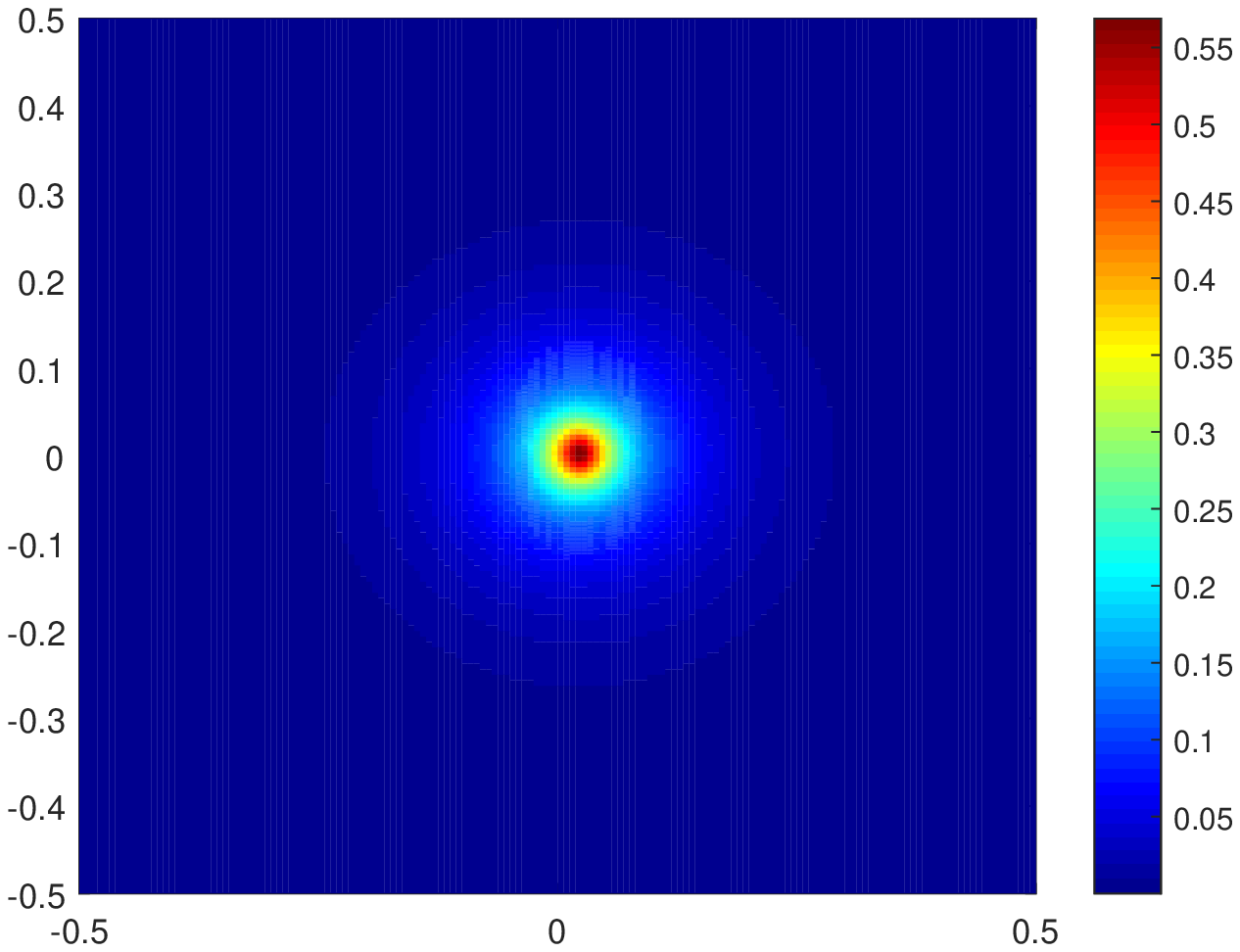}
		\includegraphics[width=0.45\textwidth,height=4.5cm]{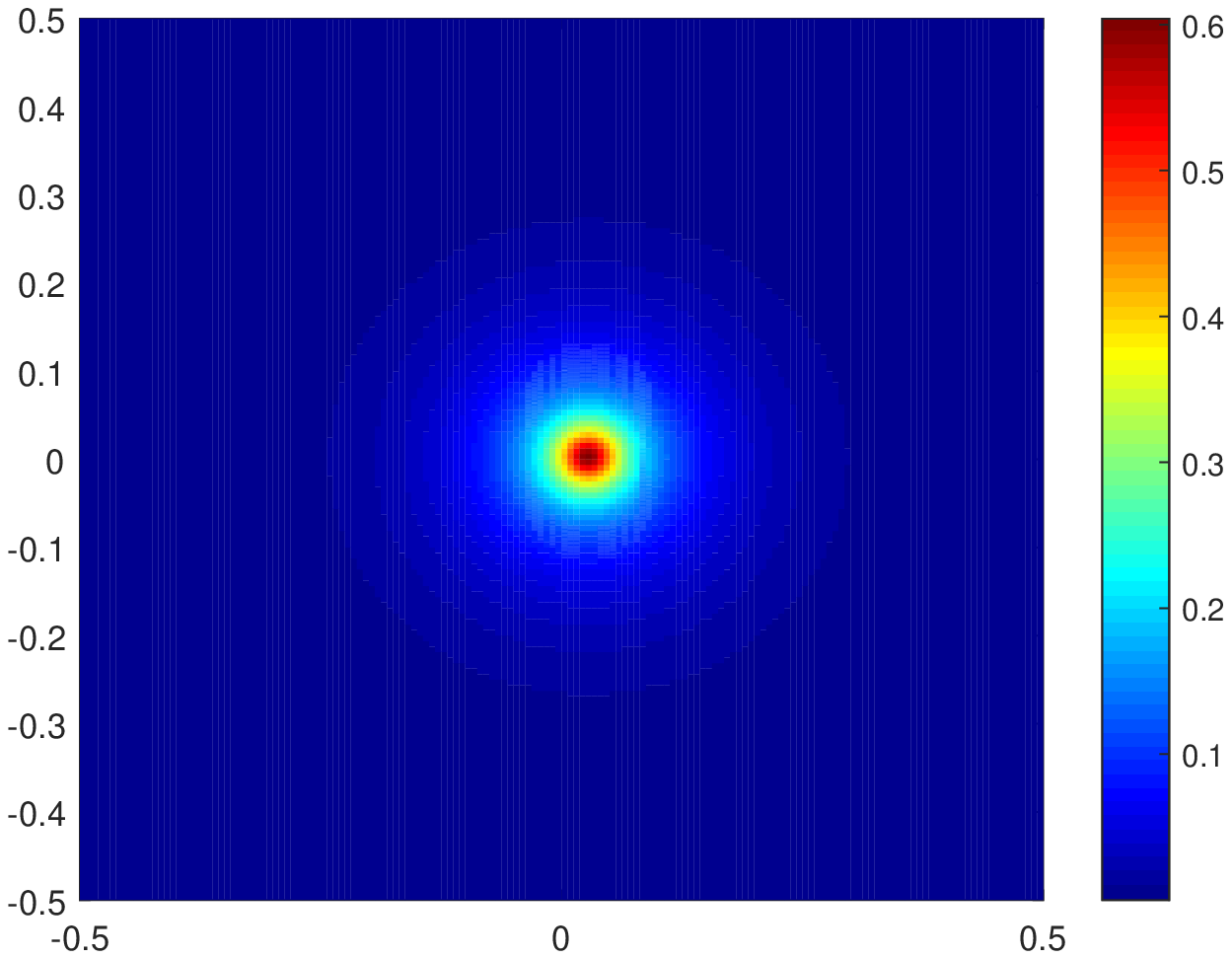}
		\caption{The case of $64$ holes: the scattered field is computed in the plane $z=0.025$. The left one shows the numerical solution to \eqref{ibvp}, and the right one shows the numerical solution to \eqref{MW-1}.}
		\label{many64_2}
	\end{center}
\end{figure}

In our setting, the scaled capacitance is $\overline C=4\pi$. Using the finite element method, we solve the original scattering problem \eqref{ibvp} and the effective problem \eqref{MW-1}, and then compare their solutions in the planar domain $\{(x,\,y,\,z): -0.5\leq x,\,y \leq 0.5,\,z=0.025\}$; see Figure \ref{many64_2}. The numerical result greatly support our theoretical result in Theorem \ref{Main2}.

\bigskip
{\bf Acknowledgement:} This work is supported by National Natural Science Foundation of China (No. 11671082) and the Austrian Science Fund(FWF): P28971-N32. The authors thank Mr. Yi Li for his contributions on the numerics. The second author would also like to thank Radon institute (RICAM), Austrian Academy of Sciences, for the friendly atmosphere during his visits.


\begin{thebibliography}{99}
\bibitem{A-C-K-S-15} B. Ahmad, D.P. Challa, M. Kirane, M. Sini, The equivalent refraction index for the acoustic scattering by many small holes: with error estimates, J. Math. Anal. Appl., 424 (2015), 563--583.

\bibitem{A-C-C-S-18} H. Ammari, D.P. Challa, A.P. Choudhury, M. Sini, The equivalent media generated by bubbles of high contrasts: Volumetric metamaterials and metasurfaces, arXiv: 1811.02912 (to appear in Multiscale Model. Simul.).

\bibitem{A-C-C-S-19} H. Ammari, D.P. Challa, A.P. Choudhury, M. Sini, The point-interaction approximation for the fields generated by contrasted bubbles atarbitrary fixed frequencies, J. Differential Equations, 267 (2019), 2104--2191.

\bibitem{A-G-H2007} H. Ammari, R. Griesmaier, M. Hanke, Identification of small inhomogeneities: asymptotic factorization, Math. Comput., 76 (2007), 1425--1448.

\bibitem{A-K:2007} H. Ammari, H. Kang, Polarization and Moment Tensors, With Applications to Inverse Problems and Effective Medium Theory, Springer, New York, 2007.

\bibitem{B-H1986_1} A. Bamberger, T. Ha Duong, Formulation variationnelle espace-temps pour le calcul par potentiel retard\'e de la diffraction d'une onde acoustique. I, Math. Methods Appl. Sci., 8 (1986), 405--435.

\bibitem{Bendali-et.al:2015} A. Bendali, P.H. Cocquet, S. Tordeux, Approximation by multipoles of the multiple acoustic scattering by small obstacles in three dimensions and application to the Foldy theory of isotropic scattering, Arch. Rational Mech. Anal., 219 (2016), 1017--1059.

\bibitem{B-L-P:1978} A. Bensoussan, J.-L. Lions, G. Papanicolaou, Asymptotic Analysis for Periodic Structures, North-Holland Publishing Co., Amsterdam, 1978.

\bibitem{BS2018} A. Bouzekri, M. Sini, The Foldy-Lax approximation for the full electromagnetic scattering by small conductive bodies of arbitrary shapes, Multiscale Model. Simul., 17 (2019), 344--398.


\bibitem{C-D2016} B. Cabarrubias, P. Donato, Homogenization of some evolution problems in domains with small holes, Electron. J. Differential Equations, 2016, Paper No. 169, 26 pp.


\bibitem{C-M-S-17} D.P. Challa, A. Mantile, M. Sini, Characterization of the equivalent acoustic scattering for a cluster of an extremely large number of small holes, arXiv: 1711.05003v1 (to appear in Asymptotic Analysis).

\bibitem{C-S2014} D.P. Challa, M. Sini, On the justification of the Foldy-Lax approximation for the acoustic scattering by small rigid bodies of arbitrary shapes, Multiscale Model. Simul., 12 (2014), 55--108.

\bibitem{Costabel2003} M. Costabel, Time-dependent problems with the boundary integral equation method, Encyclopedia of Computational Mechanics, Erwin Stein, Renee de Borst and Thomas Hughes, eds., John Wiley, New York, 2003.

\bibitem{D-G2008}  P. Donato, F. Gaveau, Homogenization and correctors for the wave equation in non periodic perforated domains, Netw. Heterog. Media, 3 (2008), 97--124. 

\bibitem{D-Z2012} P. Donato, Z. Yang, The periodic unfolding method for the wave equation in domains with holes,
Adv. Math. Sci. Appl., 22 (2012), 521--551.

\bibitem{G-R1986} V. Girault, P.-A. Raviart, Finite Element Methods for Navier-Stokes Equations: Theory and Algorithms, Springer-Verlag, Berlin, 1986.

\bibitem{Ha2003} T. Ha-Duong, On retarded potential boundary integral equations and their discretisations, Topics in computational wave propagation, 301--336, Lect. Notes Comput. Sci. Eng., 31, Springer, Berlin, 2003.

\bibitem{J-K-O:1994} V.V. Jikov, S.~M. Kozlov, O.~A. Ole{\u\i}nik, Homogenization of Differential Operators and Integral
Functionals, Springer-Verlag, Berlin, 1994.

\bibitem{K-P2017} D.V. Korikov, B.A. Plamenevskii, Asymptotics of solutions of the stationary and nonstationary Maxwell systems in a domain with small cavities, St. Petersburg Math. J., 28 (2017), 507--554.

\bibitem{L-M2015} A. Lechleiter, P. Monk, The time-domain Lippmann-Schwinger equation and convolution quadrature, Numerical Methods for Partial Differential Equations, 31 (2015), 517--540.

\bibitem{Lubich1994} Ch. Lubich, On the multistep time discretization of linear initial-boundary value problems and their boundary integral equations, Numer. Math., 67 (1994), 365--389.

\bibitem{M-K:2006} V.A. Marchenko, E.Y. Khruslov, Homogenization of Partial Differential Equations, Birkh{\"a}user Boston Inc., Boston, MA, 2006.

\bibitem{Martin:2006} P.A. Martin, Multiple Scattering. Interaction of Time-harmonic Waves with $N$ Obstacles, Cambridge University Press, Cambridge, 2006.

\bibitem{M-M-N2017} V.G. Maz'ya, A.B. Movchan, M.J. Nieves, Eigenvalue problem in a solid with many inclusions: asymptotic analysis, Multiscale Model. Simul., 15 (2017), 1003--1047.

\bibitem{M-N-P2000} V.G. Maz'ya, S. Nazarov, B. Plamenevskij, Asymptotic Theory of Elliptic Boundary Value Problems in Singularly Perturbed Domains. I, II, Birkh\"auser, Basel, 2000.

\bibitem{Nie2017} M.J. Nieves, Asymptotic analysis of solutions to transmission problems in solids with many inclusions, SIAM J. Appl. Math., 77 (2017), 1417--1443.

\bibitem{Ramm2015} A.G. Ramm, Scattering of electromagnetic waves by many small perfectly conducting or impedance bodies, J. Math. Phys., 56 (2015), 091901.

\bibitem{Sayas2016} F.-J. Sayas, Retarded Potentials and Time Domain Boundary Integral Equations, Springer, New York, 2016.

\bibitem{S-W2019} M. Sini, H. Wang, Estimation of the heat conducted by a cluster of small cavities and characterization of the equivalent heat conduction, Multiscale Model. Simul., 17 (2019), 1214--1251.

\end{thebibliography}
\end{document}